\documentclass[a4paper]{amsart}

\usepackage[english]{babel}
\usepackage{array}
\usepackage[utf8]{inputenc}
\usepackage{graphicx}
\usepackage{amsmath,latexsym,amsfonts,amssymb,amscd,euscript,epic,eepic,times,amsthm,enumerate}
\usepackage{xypic}
\usepackage{nicefrac}
\usepackage{faktor}
\usepackage{pdfsync}
\usepackage{hyperref}
\usepackage{url}
\usepackage[all]{xy}

\usepackage{hyperref}
\usepackage{lscape}
\usepackage{array}

\newtheorem{theorem}{Theorem}[section]
\newtheorem{lemma}[theorem]{Lemma}
\newtheorem{proposition}[theorem]{Proposition}
\newtheorem{corollary}[theorem]{Corollary}
\newtheorem{definition}[theorem]{Definition}
\newtheorem{remark}[theorem]{Remark}

\newcommand{\bpr}{\begin{proposition}}
\newcommand{\epr}{\end{proposition}}
\newcommand{\btm}{\begin{theorem}}
\newcommand{\etm}{\end{theorem}}
\newcommand{\bco}{\begin{corollary}}
\newcommand{\eco}{\end{corollary}}
\newcommand{\blm}{\begin{lemma}}
\newcommand{\elm}{\end{lemma}}
\newcommand{\bdf}{\begin{definition}}
\newcommand{\edf}{\end{definition}}
\newcommand{\bpm}{\begin{pmatrix}}
\newcommand{\epm}{\end{pmatrix}}
\newcommand{\beq}{\begin{equation}}
\newcommand{\eeq}{\end{equation}}
\newcommand{\bit}{\begin{itemize}}
\newcommand{\eit}{\end{itemize}}
\newcommand{\brm}{\begin{remark}}
\newcommand{\erm}{\end{remark}}
\newcommand{\bpf}{\begin{proof}}
\newcommand{\epf}{\end{proof}}

\newcommand{\lie}{\mathfrak}
\newcommand{\nr}{\textnormal}
\newcommand{\wt}{\widetilde}

\newcommand{\ol}{\overline}

\newcommand{\CC}{\mathbb{C}}

\newcommand{\RR}{\mathbb{R}}

\newcommand{\ZZ}{\mathbb{Z}}

\newcommand{\Aa}{\mathcal{A}}
\newcommand{\Bb}{\mathcal{B}}
\newcommand{\Cc}{\mathcal{C}}

\newcommand{\Ee}{\mathcal{E}}

\newcommand{\Hh}{\mathcal{H}}

\newcommand{\Mm}{\mathcal{M}}

\newcommand{\Oo}{\mathcal{O}}
\newcommand{\Pp}{\mathcal{P}}

\newcommand{\Rr}{\mathcal{R}}

\newcommand{\Tt}{\mathcal{T}}

\newcommand{\gG}{\lie{g}}
\newcommand{\hH}{\lie{h}}
\newcommand{\kK}{\lie{k}}

\newcommand{\lL}{\lie{l}}
\newcommand{\pP}{\lie{p}}

\newcommand{\sS}{\lie{s}}

\newcommand{\zZ}{\lie{z}}

\newcommand{\quotient}[2]{{\raisebox{.2em}{\thinspace $#1$}\left / \raisebox{-.15em}{ $#2$}\right.}}
\newcommand{\git}[2]{{\raisebox{.2em}{\thinspace $#1$}\left /\!\!/ \raisebox{-.15em}{$#2$}\right.}}

\newcommand\Quotient[2]{
        \mathchoice
            {
                \text{\raise1ex\hbox{\thinspace #1}\Big{/} \lower1ex\hbox{#2} \thinspace}%
            }
            {
                #1\,/\,#2
            }
            {
                #1\,/\,#2
            }
            {
                #1\,/\,#2
            }
    }

\newcommand\GIT[2]{
        \mathchoice
            {
                \text{\raise1ex\hbox{\thinspace #1}\Big{/}\!\!\!\!\Big{/} \lower1ex\hbox{#2} \thinspace}%
            }
            {
                #1\,/\,#2
            }
            {
                #1\,/\,#2
            }
            {
                #1\,/\,#2
            }
    }

\newcommand{\morph}[6]{\begin{array}{cccc} #1 \quad : &  #2  & \stackrel{#6}{\lra} &  #3  \\  & #4 &\longmapsto & #5  \end{array}}

\newcommand{\map}[5]{\begin{array}{cccc}   #1  & \stackrel{#5}{\lra} &  #2  \\  #3 & \longmapsto & #4  \end{array}}

\newcommand{\injection}[4]{\begin{array}{cccc}   #1  & \hookrightarrow &  #2  \\  #3 &\longmapsto & #4  \end{array}}

\newcommand{\mapp}[7]{\begin{array}{cccccc}   #1  & \lra &  #2 & \stackrel{#7}{\lra} &  #3   \\  #4 &\longmapsto & #5 & \longmapsto & #6  \end{array}}

\newcommand{\momorph}[9]{\begin{array}{cccc} #1 \quad : &  #2  & \stackrel{#8}{\lra} &  #3  \\  &  #4  & \stackrel{#9}{\lra} &  #5  \\ & #6 &\longmapsto & #7  \end{array}}

\newcommand{\dif}{\thinspace d}
\newcommand{\qua}{\thinspace}
\newcommand{\dolbeault}{\overline{\partial}}
\newcommand{\Slash}{/\!\!/}

\newcommand{\lra}{\longrightarrow}

\DeclareMathOperator{\id}{id}

\DeclareMathOperator{\ad}{ad}

\DeclareMathOperator{\Sym}{Sym}

\DeclareMathOperator{\Pic}{Pic}

\DeclareMathOperator{\im}{im}

\DeclareMathOperator{\codim}{codim}

\DeclareMathOperator{\GL}{GL}
\DeclareMathOperator{\PGL}{PGL}
\DeclareMathOperator{\SL}{SL}

\DeclareMathOperator{\U}{U}
\DeclareMathOperator{\SU}{SU}

\DeclareMathOperator{\Hom}{Hom}

\DeclareMathOperator{\aut}{aut}

\title{Higgs bundles over elliptic curves for complex reductive groups}

\author{Emilio Franco}
\address{Emilio Franco \\ CMUP (Centro de Matem\'atica da Universidade do Porto) \\ Universidade do Porto \\ Rua do Campo Alegre 1021/1055 \\ 4169-007, Porto (Portugal)}
\email{emilio.franco@fc.up.pt}

\author{Oscar Garcia-Prada}
\address{Oscar Garcia-Prada \\ ICMAT (Instituto de Ciencias Matem\'aticas) \\ CSIC-UAM-UC3M-UCM \\ 
Calle Nicol\'as Cabrera 15 \\ 28049 Madrid (Spain)}
\email{oscar.garcia-prada@icmat.es}

\author{P. E. Newstead}
\address{P. E. Newstead \\ Department of Mathematical Sciences \\ University of Liverpool \\ 
Peach Street \\ Liverpool L69 7ZL (United Kingdom)}
\email{newstead@liverpool.ac.uk}

\begin{document}

\date{\today}

\keywords{Higgs bundles, elliptic curves, Hitchin map.}

\subjclass[2010]{14H60, 14D20, 14H52} 

\thanks{First author partially supported by Consejo Superior de Investigaciones Cient\'ificas (CSIC) through JAE-Predoc grant program, German Research Foundation through the project SFB 647 and Fundação de Amparo à Pesquisa do Estado de São Paulo (FAPESP)  through grant 2012/16356-6. First and second authors partially supported by the Ministerio de Econom\'ia y Competitividad of Spain through Project MTM2010-17717 and Severo Ochoa Excellence Grant. The three authors thank the Isaac Newton Institute in Cambridge --- which they visited while preparing the paper --- for the excellent conditions provided.}

\begin{abstract}
We study Higgs bundles over an elliptic curve with complex reductive structure group, describing the (normalization of) its moduli spaces and the associated Hitchin fibration. The case of trivial degree is covered by the work of Thaddeus in 2001. Our arguments are different from those of Thaddeus and cover arbitrary degree.
\end{abstract}

\maketitle

\tableofcontents

\section{Introduction}
\label{sc introduction}

An {\it elliptic curve} is a pair $(X,x_0)$ where $X$ is a smooth complex projective curve of genus $1$ and $x_0$ is a distinguished point on it. By abuse of notation, we usually refer to an elliptic curve simply as $X$. Let $G$ be a connected complex reductive Lie group. A {\it $G$-Higgs bundle} over $X$ is a pair $(E,\Phi)$ where $E$ is a principal $G$-bundle over $X$ and $\Phi$, called the {\it Higgs field}, is a section of the adjoint bundle twisted by the canonical bundle of the curve. The canonical bundle of an elliptic curve is trivial, $\Omega^1_X \cong \Oo_X$, so $\Phi \in H^0(X,E(\gG))$. These objects were defined by Hitchin \cite{hitchin-self_duality_equations} over a smooth projective curve of any genus and the existence of their moduli spaces $\Mm(G)_d$ (here $d\in\pi_1(G)$ is a topological invariant known as the degree) follows from Simpson \cite{simpson1, simpson2} (the existence of $\Mm(\SL(2,\CC))$ was first given in \cite{hitchin-self_duality_equations} and the case of $\GL(n,\CC)$ was also given by Nitsure \cite{nitsure}).

A major result of the theory of $G$-Higgs bundles is the non-abelian Hodge correspondence which was proved by Hitchin \cite{hitchin-self_duality_equations}, Do\-naldson \cite{donaldson}, Simpson \cite{simpson-hb&ls, simpson1, simpson2} and Corlette \cite{corlette}. It is a generalization of the Narasimhan--Seshadri--Ramanathan Theorem \cite{narasimhan&seshadri, ramanathan_stable} to the non-unitary case and states the existence of a chain of homeomorphisms between the moduli space of $G$-Higgs bundles, the moduli space of $G$-bundles with projectively flat connections $\Cc(G)_d$ and the moduli space of representations $\Rr(G)_d$ of the curve,
\[
\Mm(G)_d \stackrel{homeo}{\cong} \Cc(G)_d \stackrel{homeo}{\cong} \Rr(G)_d.
\]

The Hitchin fibration was defined by Hitchin \cite{hitchin-duke} using a basis $p_1, \dots, p_\ell$ of the invariant polynomials of the Lie algebra $\gG$,
\[
\map{\Mm(G)_d}{B_G \cong \bigoplus H^0(X,(\Omega^1_X)^{\otimes \deg(p_i)})}{(E,\Phi)}{(p_1(\Phi), \dots, p_\ell(\Phi)).}{ }
\]
A more canonical definition of the Hitchin fibration was provided by Donagi \cite{donagi} re\-de\-fi\-ning the Hitchin base $B_G$ as the space of cameral covers $H^0(X, (\gG \otimes \Omega^1_X)\Slash G)$. Another ground-breaking result of the theory of Higgs bundles says that, under this fibration, the space of $G$-Higgs bundles is an algebraically completely integrable system \cite{hitchin-duke, faltings, donagi}.

\

In 1957 Atiyah \cite{atiyah} studied vector bundles over an elliptic curve $X$ leading to an identification of the moduli space of vector bundles $M(\GL(n,\CC))_d$ with $\Sym^h X$, where $h$ is the greatest common divisor of $n$ and $d$. Some forty years later, Laszlo \cite{laszlo} and Friedman, Morgan and Witten \cite{friedman&morgan, friedman&morgan&witten}, gave a description of the moduli space of $G$-bundles $M(G)_d$ (\cite{laszlo} only deals with $M(G)_0$) as the quotient 
\beq
\label{eq description of MG}
M(G)_d \cong \quotient{(X \otimes_\ZZ \Lambda_{G,d})}{W_{G,d}}
\eeq
where $\Lambda_{G,d}$ is a certain lattice, $W_{G,d}$ is a finite group acting on $\Lambda_{G,d}$ and $X \otimes_\ZZ \Lambda_{G,d}$ is the tensor product over $\ZZ$ (recall that $X$ is an abelian variety and therefore has a natural $\ZZ$-module structure). When $G$ is simply connected (and therefore $d = 0$), $\Lambda_{G,0} = \Lambda$ is the coroot lattice and $W_{G,0} = W$ is the Weyl group of $G$. In this case, by a result of Looijenga \cite{looijenga} (see also \cite{bernstein&shavartzman}), $M(G)_0$ is isomorphic to a weighted projective space. This isomorphism was obtained directly by Friedman and Morgan \cite{friedman&morgan_2} working with deformations of unstable $G$-bundles (see also \cite{helmke&slowody}). 

The construction of the isomorphism (\ref{eq description of MG}) relies on two facts. The first one is the des\-crip\-tion of the moduli space of unitary representations $R(G)_d$ achieved by Schweigert \cite{schweigert} and more generally by Borel, Friedman and Morgan \cite{borel&friedman&morgan}. By the Narasimhan--Seshadri--Ramanathan Theorem, $R(G)_d$ is homeomorphic to $M(G)_d$. This shows that an appropiate morphism from $(X \otimes_\ZZ \Lambda_{G,d})/W_{G,d}$ to $M(G)_d$ is bijective. The other key result is the fact that $M(G)_d$ is a normal projective variety, which allows us to apply Zariski's Main Theorem, proving that the previous bijective morphism is indeed an isomorphism. 

In this paper, we describe $\Mm(G)_d$ for any complex reductive group $G$, thus generalising \cite{franco&oscar&newstead}, where the authors studied these objects when $G$ is a classical group.

\

The results of this paper are structured as follows. After reviewing in Section \ref{sc G-bundles} the theory of unitary representations and $G$-bundles over an elliptic curve, we prove in Section \ref{sc chapter G-Higgs bundles} that a $G$-Higgs bundle is (semi)stable if and only if the underlying $G$-bundle is (semi)stable [Propositions \ref{pr E Phi semistable iff E semistable for G-Higgs bundles} and \ref{pr E Phi stable iff E stable for G-Higgs bundles}]. This fact shows the existence of a projection [Corollary \ref{co MmG onto MG}]
\beq
\label{eq MmG to MG}
\Mm(G)_d \lra M(G)_d
\eeq
and, combined with the results of \cite{borel&friedman&morgan}, implies that every polystable $G$-Higgs bundle of degree $d$ reduces to a unique (up to conjugation) Jordan--Hölder Levi subgroup $L_{G,d}$ [Proposition \ref{pr unique Jordan-Holder Levi subgroup for G-Higgs bundles}]. This allows us to give a complete description of the polystable $G$-Higgs bundles [Corollaries \ref{co description polystable G-Higgs bundles} and \ref{co automorphisms groupn of polystable G-Higgs bundles}]. Using this description, we construct a family $\Hh_{G,d}$ of polystable $G$-Higgs bundles of degree $d$ parametrized by $T^*X \otimes_\ZZ \Lambda_{G,d}$. Every polystable $G$-Higgs bundle can be constructed starting from a Higgs bundle for an abelian group [Remark \ref{rm abelianization}], which shows that the non-abelian Hodge correspondence is not entirely non-abelian in the elliptic case. Next, we show that the morphism associated to the family $\Hh_{G,d}$ factors through a bijective morphism and, using Zariski's Main Theorem, this gives us a description of the normalization of the moduli space [Theorem \ref{tm MmG is the normalization of MmMG}]
\beq
\label{eq normalization of Mm}
\ol{\Mm(G)}_d \cong \quotient{(T^*X \otimes_\ZZ \Lambda_{G,d})}{W_{G,d}}.
\eeq

It is not known whether $\Mm(G)_d$ is a normal quasiprojective variety (see \cite[Section 3.4]{franco&oscar&newstead} for a detailed discussion), so we can not apply the method used to prove (\ref{eq description of MG}) since the hypothesis of Zariski's Main Theorem requires the normality of the target. By means of this bijection and the quotient (\ref{eq MmG to MG}), we define a natural orbifold structure on $M(G)_d$ and the projection (\ref{eq MmG to MG}) corresponds with the projection of the associated cotangent orbifold bundle [Remark \ref{rm orbifold structure}].

In Section \ref{sc the hitchin fibration}, we study the Hitchin fibration and we obtain that it corresponds to the projection [Proposition \ref{pr description of the Hitchin fibration}]
\[
\quotient{(T^*X \otimes_\ZZ \Lambda_{G,d})}{W_{G,d}} \lra \quotient{(\CC \otimes_\ZZ \Lambda_{G,d})}{W_{G,d}}
\]
induced by the obvious projection from $T^*X \cong X \times \CC$  to $\CC$.  This gives us an explicit description of (the normalization of) all the fibres of the Hitchin fibration and, more concretely, the generic ones [Corollary \ref{co generic Hitchin fibre}].

In Section \ref{sc representations}, we use the non-abelian Hodge correspondence and our description of $G$-Higgs bundles to extend the results of \cite{borel&friedman&morgan} about unitary representations of surface groups of an elliptic curve to reductive representations of this surface group [Corollaries \ref{co description of rank zero c-pairs for G} and \ref{co B xy = B_c and D_c has a unique rk0 c-pair for G}]. This allows us to construct a bijective morphism to the moduli space $\Rr(G)_d$ of representations and then the normalization of the moduli space is [Corollary \ref{co olzeta for G}]
\beq
\label{eq normalizarion of Rr}
\ol{\Rr(G)}_d \cong \quotient{(\CC^* \times \CC^*) \otimes_\ZZ \Lambda_{G,d}}{W_{G,d}}.
\eeq

In Section \ref{sc connections}, we study the moduli space $\Cc(G)_{d}$ of $G$-bundles with projectively flat connections. Using only the Narasimhan--Seshadri--Ramanathan Theorem and the fact that the underlying $G$-bundle of a polystable $G$-Higgs bundle is also polystable, we observe a splitting of the Hitchin equations [Proposition \ref{co splitting of the hitchin equations}] which simplifies the proof of the Hitchin--Kobayashi correspondence over elliptic curves [Corollary \ref{co Hitchin Simpson correspondence}, Remark \ref{rm Hitchin-Simpson easier}]. We obtain that the normalization of the moduli space is [Theorem \ref{tm DdG bijective Xsharp otimes Lambda}]
\beq
\label{eq normalization of Cc}
\ol{\Cc(G)}_d \cong \quotient{(X^\sharp \otimes_\ZZ \Lambda_{G,d})}{W_{G,d}},
\eeq
where we recall that $X^\sharp$ is the moduli space of rank $1$ local systems on $X$.

\

In the trivial degree case, (\ref{eq normalization of Mm}), (\ref{eq normalizarion of Rr}) and (\ref{eq normalization of Cc}) become
\[
\ol{\Mm(G)}_0 \cong \quotient{(T^*X \otimes_\ZZ \Lambda)}{W},
\]
\[
\ol{\Rr(G)}_0 \cong \quotient{((\CC^* \times \CC^*) \otimes_\ZZ \Lambda)}{W}
\]
and 
\[
\ol{\Cc(G)}_0 \cong \quotient{(X^\sharp \otimes_\ZZ \Lambda)}{W},
\]
where $W$ is the Weyl group of $G$ and $\Lambda$ is the lattice given by the kernel of the exponential restricted to the Cartan subalgebra (i.e. the fundamental group of the Cartan subgroup). This was obtained by Thaddeus \cite{thaddeus} in 2001. Our arguments are different from those of Thaddeus and work for arbitrary $d$.

When $G = \GL(n,\CC)$ or $\SL(n,\CC)$ (for any $n$, not only for $n\leq 4$ as stated in \cite{franco&oscar&newstead}) one actually obtains an isomorphism since the target is normal. In these cases, 
\[
\Rr(G)_0 := \git{\Hom(\ZZ \oplus \ZZ, G)}{G} \subset \git{\{ x,y \in \lie{g} : [x,y]= 0 \}}{G}
\]
is normal due to \cite[Section 0.2]{joseph} (although the hypothesis of \cite{joseph} requires $G$ to be semisimple, the proof can be extended to $\GL(n,\CC)$ as in \cite[Corollary 7.4]{levasseur}). Normality of $\Mm(G)_0$ and $\Cc(G)_0$ follow from the Isosingularity Theorem \cite[Theorem 10.6]{simpson2} and normality of $\Rr(G)_0$. The corresponding results for $\Rr(G)_0$ for general reductive $G$ constitute a long-standing open problem and the case of $\Rr(G)_d$ is still more uncertain. Indeed it is not even clear whether the moduli spaces are reduced.

\

We work in the category of algebraic schemes over $\CC$. Unless otherwise stated, all the bundles considered are algebraic bundles.

\

\noindent {\it Acknowledgements.} This article is a modified version of the third part of the PhD thesis of the first author \cite{phD}, prepared under the supervision of the second and third authors at ICMAT (Madrid). The first author wishes to thank the second and third authors for their teaching, help and encouragement. 

The three authors thank the anonymous referees of a previous version of the paper for valuable suggestions. In particular, the referees drew our attention to the reference \cite{thaddeus} and also to \cite{joseph} and \cite{levasseur}. Thanks are also due to the referees of the current version as the corrections they suggested have improved the paper.

\section{Review on $G$-bundles and unitary representations over elliptic curves}
\label{sc G-bundles}

\subsection{Review on the abelian case}
\label{sc the abelian case}

If $X$ is an elliptic curve, the Abel--Jacobi map gives an isomorphism $X \cong \Pic^1(X)$. Fixing a point $x_0 \in X$ and tensoring by $\Oo(x_0)^{-1}$ one obtains $\varsigma_{1,0}: X \stackrel{\cong}{\lra} \Pic^0(X)$, which induces an abelian group structure on $X$. There is a unique Poincaré bundle $\Pp \to X \times \Pic^0(X)$ such that its restriction to the slice $\{ x_0 \} \times \Pic^0(X)$ is the trivial line bundle.

Let $S$ be a compact connected abelian group and let $S^\CC$ be its complexification. The universal cover of $S$ (resp. $S^\CC$) is its Lie algebra $\sS$ (resp. $\sS^\CC$) and the covering map is the exponential $\exp : \sS \to S$ (resp. $\sS^\CC \to S^\CC$). By construction, the kernels of the two maps coincide and we write 
\[
\Lambda_S := \Lambda_{S^\CC} := \ker \exp,
\]
which is a lattice in $\sS \subset \sS^\CC$. Note that the fundamental groups $\pi_1(S)$ and  $\pi_1(S^\CC)$ coincide since both are identified with the kernel of the exponential map.

Every element $\gamma \in \Lambda_{S}$ defines a cocharacter $\theta : \CC^* \to S^\CC$ that restricts to $\theta : \U(1) \to S$. Let $\Bb = \{ \gamma_1, \dots, \gamma_k \}$ be a basis of $\Lambda_{S}$ and let $\{ \theta_1, \dots, \theta_k \}$ be the associated cocharacters. These give  isomorphisms
\beq
\label{eq definition of Theta}
\momorph{\Theta_S}{\CC^* \otimes_\ZZ \Lambda_S}{S^\CC}{\U(1) \otimes_\ZZ \Lambda_S}{S}{\sum_{i=1}^k \ell_i \otimes_\ZZ \gamma_i}{\Pi_{i=1}^k \theta_{i}(\ell_i),}{\cong}{\cong}
\eeq
where $\ell_i \in \CC^*$ (resp. $\U(1)$), and
\beq
\label{eq definition of dTheta}
\momorph{d \Theta_S}{\CC \otimes_\ZZ \Lambda_S}{\sS^\CC}{\RR \otimes_\ZZ \Lambda_S}{\sS}{\sum_{i=1}^k (\lambda_i \otimes_\ZZ \gamma_i)}{\sum_{i=1}^k \lambda_i \cdot \gamma_i,}{\cong}{\cong}
\eeq
where $\lambda_i \in \CC$ (resp. $\RR$).

Using (\ref{eq definition of Theta}) and fibre products of the Poincaré bundle $(\id \times \varsigma_{1,0})^*\Pp \to X \times X$, one can construct a family of $S^\CC$-bundles with trivial degree,
\beq
\label{eq poincare bundle}
\Pp_{S} \lra X \times (X \otimes_\ZZ \Lambda_{S}),
\eeq
whose restriction to the slice $\{ x_0 \} \times (X \otimes_\ZZ \Lambda_{S})$ is the trivial $S^\CC$-bundle over $(X \otimes_\ZZ \Lambda_{S})$. 

Among other references, the following result is contained in \cite[Theorem 9.6]{simpson2} (recall that for an elliptic curve $X \cong \Pic^0(X)$). 

\btm
\label{tm MG for G abelian}
Let $S^\CC$ be an abelian, connected complex Lie group. Then, the moduli space of topologically trivial $S^\CC$-bundles over the elliptic curve $X$ is
\[
M(S^\CC)_0 \cong X \otimes_\ZZ \Lambda_{S}.
\]
\etm

\subsection{Notation and some results on Lie groups}
\label{sc notation and results on Lie groups}

We refer to \cite{phD} for an expanded version of this section. Let $G$ denote a compact (resp. complex reductive) connected Lie group. We set some notation:
\begin{itemize}
\item $Z_0$ denotes the connected component of the identity of the center $Z_G(G)$ of the group,
\item $p : D \to [G,G]$ denotes the universal covering of the semisimple group $[G,G]$,
\item $F := Z_0 \cap [G,G]$,
\item $C := p^{-1}(F) \subset Z_D(D)$,
\item $\tau : C \to Z_0$ denotes the homomorphism given by the inclusion $F \hookrightarrow Z_0$.
\item $\ol{G} := G / F$,
\item $\ol{Z} := Z_0 / F$,
\item $\ol{D} := D / C$ or equivalently $[G,G]/F$,
\item $H \subset G$ denotes a maximal torus (resp. Cartan subgroup) with Lie algebra $\hH$,
\item $H' \subset D$ denotes a maximal torus (resp. Cartan subgroup) with Lie algebra $\hH' = [\hH,\hH]$,
\item $W = N_G(H) / Z_G(H) = N_D(H') / Z_D(H')$ denotes the Weyl group.
\end{itemize}

Note that we have natural isomorphisms
\beq
\label{eq G cong Z_0 times_tau D}
G \cong Z_0 \times_\tau D
\eeq
and 
\beq
\label{eq olG = olZ times olD}
\ol{G} \cong \ol{Z} \times \ol{D}.
\eeq
The finite covering $G \to \ol{G}$ induces an injection
\beq
\label{eq definition of q^pi}
\injection{\pi_1(G)}{\pi_1(\ol{Z}) \times \pi_1(\ol{D})}{d}{(u,c).}
\eeq
Since $D$ is simply connected and $C$ finite, we have
\[
\pi_1(\ol{D}) = C.
\]

Let us suppose for simplicity that $D$ is a simple compact Lie group (resp. simple complex Lie group). Take an alcove $A \subset \hH'$ containing the origin. For $c \in Z_D(D)$, we know (see for instance \cite{broker&tomDieck}) that there is a vertex $a_c$ of the alcove $A$  such that $c = \exp(a_c)$. We see that $A - a_c$ is another alcove contaning the origin. Hence there is a unique element $\omega_c \in W$ such that
\[
A - a_c = \omega_c(A).
\]
In the trivial case we obviously have $\omega_0 = \id$.

We denote the connected component of the fixed point set of the action of $\omega_c$ on $H$ by
\beq
\label{eq definition of S_c}
S_c:= (H^{\omega_c})_0.
\eeq
Let us take its normalizer $N_G(S_c)$ and define the quotient 
\beq
\label{eq definition of W_c}
W_c := \quotient{N_G(S_c)}{Z_G(S_c)} = \quotient{N_G(\hH^{\omega_c})}{Z_G(\hH^{\omega_c})}.
\eeq
When $c$ is the identity, one recovers the usual Weyl group $W$.

We define 
\beq
\label{eq definition of L_c}
L_c:= Z_G(S_c). 
\eeq
Since $L_c$ is the centralizer of a torus, we know that it is connected. One can easily check that $N_G(S_c) = N_G(L_c)$ and therefore $W_c = \quotient{N_G(L_c)}{L_c}$.

By \cite[Lemma 2.1.1]{borel&friedman&morgan} and \cite[Proposition 3.4.4]{borel&friedman&morgan}, $D_c = [L_c,L_c]$ is simply connected. Define $F_c = S_c \cap D_c$ and note that $S_c$ is the centre of $L_c$. By (\ref{eq G cong Z_0 times_tau D}), we have $L_c \cong S_c \times_{F_c} D_c$. Note, by (\ref{eq definition of q^pi}), that $\pi_1(L_c)$ injects into $\pi_1(\ol{S_c}) \times \pi_1(D_c/F_c)$, where 
\beq
\label{eq definition of olS_c and olD_c}
\ol{S_c}:= S_c / F_c.
\eeq

The inclusion $L_c \hookrightarrow G$ induces a morphism $\pi_1(L_c) \to \pi_1(G)$.

\blm
\label{lm pi_1 G injects into pi_1 L_S}
Let $d = (u,c) \in \pi_1(G)$ and let $L_c$ be associated to $c$. Then there is a unique $\ell_d \in \pi_1(L_c)$ that maps to $d$ and furthermore $\ell_d = (u,p(c))$.
\elm 

\proof
By construction, we have that $p(c) \in D_c = [L_c,L_c]$ and $p(c) \in S_c$, thus $p(c) \in F_c \subset Z_{D_c}(D_c)$. If $\ell \in \pi_1(L_c)$ is given by $(v,f) \in \sS \times Z_{D_c}(D_c)$ and it maps to $d$, then $f = p(c)$ and $v = u$, since $v \in \exp^{-1}(p(c)) \subset \exp^{-1}(F) \subset \zZ_{\gG}(\gG)$. The choice of $d$ fixes $(v,f)$, so its preimage $\ell \in \pi_1(L_c)$ is unique.
\qed

Recall that $W_c = N_G(\sS_c)/Z_G(\sS_c)$, where $\sS_c = \hH^{\omega_c}$ is the Lie algebra of $S_c$, and note that $W_c$ preserves $\Lambda_{S_c} \subset \sS_c$. This gives us an action of $W_c$ on $\U(1) \otimes_\ZZ \Lambda_{S_c}$ (resp. on $\CC^* \otimes_\ZZ \Lambda_{S_c}$) and this action commutes with the isomorphism $\Theta_{S_c}$ defined in (\ref{eq definition of Theta}).

In (\ref{eq definition of olS_c and olD_c}) we have defined $\ol{S}_c$ as $S_c/F_c$. We can check that $W_c$ preserves $F_c$, so the action of $W_c$ on $S_c$ gives a well defined action of $W_c$ on $\ol{S}_c$. Note that $\Lambda_{\ol{S}_c} = \exp_S^{-1}(F_c)$, so $W_c$ also preserves $\Lambda_{\ol{S}_c}$, inducing an action on $\U(1) \otimes_\ZZ \Lambda_{\ol{S}_c}$ (resp. on $\CC^* \otimes_\ZZ \Lambda_{\ol{S}_c}$). We can check that the action of $W_c$ commutes with $\Theta_{\ol{S}_c}$ too.

\subsection{Representations and $c$-pairs}
\label{sc representations and c-pairs}

In this section we present some results from \cite{borel&friedman&morgan} (see also \cite{phD}). We say that two elements of a Lie group $G$ {\it almost commute} if their commutator lies in the centre of the Lie group.  Let $c$ be an element of $C \subset Z_D(D)$. Suppose $a$ and $b$ are two almost commuting elements of the form $a=[(z_1,\delta_1)]_\tau$ and $b=[(z_2,\delta_2)]_\tau$, where $z_1, z_2 \in Z_0$ and $\delta_1, \delta_2 \in D$. We say that $(a,b)$ is a $c$-{\it pair} if $[\delta_1,\delta_2]=c$. Let $C(G)_c$ denote the subset of $G \times G$ of $c$-pairs. 

The fundamental group of an elliptic curve is $\pi_1(X) = \langle \alpha, \beta : [\alpha, \beta] = \id \rangle \cong \ZZ^2$. Take the universal central extension $\Gamma = \langle \alpha, \beta,\delta : [\alpha, \beta] = \delta, [\alpha, \delta] = \id, [\beta, \delta] = \id \rangle$ and define $\Gamma_\RR$ as $\Gamma \times_\ZZ \RR$. A representation $\rho : \Gamma_\RR \to G$ is \textit{central} if $\rho(\RR)$ is contained in $Z_G(G)$; since $\rho(\RR)$ is connected and contains the unit element, it is contained in $Z_0 = Z_G(G)_0$. From a central representation $\rho: \Gamma_\RR \to G$, one obtains a pair $(\nu,u)$, where $\nu : \Gamma \to G$ is such that $\nu = \rho|_\Gamma$ and $u \in \zZ_{\gG}(\gG)$ is given by $u = d \rho (1)$ and, thanks to the exponential map, $u$ can be viewed as an element of the fundamental group of $\ol{Z}$. Conversely, $(\nu, u)$ determines uniquely a central representation $\rho: \Gamma_\RR \to G$. We observe that $u \in \zZ_{\gG}(\gG)$ is an invariant of the conjugacy class of the representation $\rho$. We denote by $\Hom^c(\Gamma_\RR,G)_d$ the set of central representations with topological invariant $d$ and we define the {\it moduli space of such representations} as the GIT quotient by the conjugation action of the group
\[
\Rr(G)_d := \git{\Hom^c(\Gamma_\RR,G)_d}{G}.
\]

Every central representation $\nu : \Gamma \to G$ is completely determined by two elements of $G$, $a = \nu(\alpha)$ and $b = \nu(\beta)$. Since $\nu$ is central, $\nu(\delta) = [a,b]$ is contained in $Z_0$ and therefore in $F = Z_0 \cap [G,G]$. Take $a = [(z_1, \delta_1)]_\tau$ and $b = [(z_2, \delta_2)]_\tau$, and write $c = [\delta_1,\delta_2]$, where $\nu(\delta) = \tau(c)$. Then $(a,b)$ completely determines the representation $\nu : \Gamma \to G$ and is a $c$-pair. Furthermore, $c \in C \subset Z_D(D)$ is an invariant of the conjugacy class of the representation $\nu$.

\brm
\label{rm rho is determined by a b and u}
\nr{Every central representation $\rho : \Gamma_\RR \to G$ is determined by a $c$-pair $(a,b) \in C(G)_c$ and an element $u$ of $\zZ_{\gG}(\gG)$ that satisfies $\tau(c) = \exp(u)$. The pair $d = (u,c) \in \zZ_{\gG}(\gG) \times Z_D(D)$ is an invariant of the conjugacy class of $\rho$. Indeed $d$ is an element of $\pi_1(G)$ as indicated  by (\ref{eq definition of q^pi}).}  

\nr{For any $g \in G$, the representation $g \rho g^{-1}$ is determined by $(gag^{-1},gbg^{-1},u)$, where $(gag^{-1},gbg^{-1})$ is a $c$-pair.}
\erm

By Remark \ref{rm rho is determined by a b and u}, we see that $\Hom^c(\Gamma_\RR,G)_{(u,c)}$ can be identified with $C(G)_c$. As a consequence, the moduli space of representations of $\Gamma_\RR$ for an elliptic curve with invariant $d \in \pi_1(G)$ determined by $(u,c) \in \zZ_{\gG}(\gG) \times Z_D(D)$ coincides with the moduli space of $c$-pairs
\[
\Rr(G)_d \cong \git{C(G)_c}{G}. 
\]

Suppose now that $G$ is a connected complex reductive algebraic group and let $K$ be its maximal compact subgroup. A representation $\rho$ is {\it reductive} if and only if the Zariski closure of $\im \rho$ is a reductive group. It is proved in \cite{richardson} that the orbit $[\rho]_G$ is closed if and only if $\rho$ is a reductive representation. Denote by $\Hom^c(\Gamma_\RR,G)^+_d$ the set of central reductive representations, and by $C(G)^+_c$ the set of reductive $c$-pairs (those associated to reductive representations). Then
\beq
\label{eq reductive c-pairs are enougth}
\Rr(G)_d \cong \quotient{\Hom^c(\Gamma_\RR,G)^+_{d}}{G} \cong \quotient{C(G)^+_c}{G}.
\eeq
Note that, for $G$ compact, every representation of $\Gamma_\RR$ is reductive. So the moduli space of unitary representations is a categorical quotient
\[
R(G)_d := \Rr(K)_d \cong \quotient{C(K)_c}{K}.
\]

A representation $\rho$ is {\it irreducible} if the centralizer of its image, $Z_G(\rho)$, is equal to $Z_G(G)$. Analogously, we say that a $c$-pair $(a,b)$ is {\it irreducible} if the centralizer of its elements, $Z_G(a,b)$, is equal to $Z_G(G)$. 

\subsection{Review on unitary representations over elliptic curves}
\label{sc borel friedman morgan}

Following \cite{borel&friedman&morgan}, in this section we study the moduli space of central representations of $\Gamma_\RR$ into a compact Lie group $K$. Let $C = p^{-1}(F) = \pi_1(\ol{D})$ as defined at the beginning of Section \ref{sc notation and results on Lie groups} and set $c \in C$.

\bpr {\bf (\cite[Proposition 4.1.1]{borel&friedman&morgan}).}
\label{pr description of rank zero c-pairs}
Let $K$ be a simply connected compact semisimple Lie group.  Let $(a,b)$ be an irreducible $c$-pair in $K$. Then
\begin{enumerate}

\item \label{num H is a product of SUn} the group $K$ is a product of simple factors $K_i$, where each $K_i$ is isomorphic to $\SU(n_i)$ for some $n_i \geq 2$;

\item \label{num c is a product of generators of ZZ n} $c = (c_1, \dots, c_r)$, where each $c_i$ generates the centre of $K_i$;
  
\item \label{num a c-pair as above has rank zero} conversely, if $K$ is as in (\ref{num H is a product of SUn}) and $c$ as in (\ref{num c is a product of generators of ZZ n}), then there is an irreducible $c$-pair in $K$ and all $c$-pairs in $K$ are conjugate.

\end{enumerate}
\epr 

Recall that $L_c \cong S_c \times_{\tau_c} D_c$, where $S_c$, $L_c$ are defined in (\ref{eq definition of S_c}), (\ref{eq definition of L_c}) and $D_c=[L_c,L_c]$. 

\bpr {\bf (\cite[Proposition 4.2.1]{borel&friedman&morgan}).}
\label{pr B xy = B_c and D_c has a unique rk0 c-pair}
Let $K$ be a compact Lie group. Let $(a,b)$ be any $c$-pair. Any maximal torus of $Z_K(a,b)$ is conjugate in $K$ to $S_c$ , so $(a,b)$ is contained in $L_c$ after conjugation and, as a $c$-pair in $L_c$, is irreducible. 
\epr 

Now we have the ingredients to describe the moduli space of unitary representations. Fix $d \in \pi_1(G)$ determined by $(u,c) \in \pi_1(\ol{Z}) \times \pi_1(\ol{D})$ as described in (\ref{eq definition of q^pi}). Let us take $(\delta_1, \delta_2)$ to be one representative of the unique conjugation class of the irreducible $c$-pair in $D_c$. Consider the following continuous map
\beq
\label{eq zeta'}
\map{(S_c \times S_c)}{\Rr(K)_d}{(s_1,s_2)}{([s_1, \delta_1]_{\tau_c}, [s_2, \delta_2]_{\tau_c}).}{}
\eeq 
Using Proposition \ref{pr B xy = B_c and D_c has a unique rk0 c-pair}, one can check that (\ref{eq zeta'}) is surjective.

\brm
\label{rm rho conjugated to rho_J rho}
\nr{By Proposition \ref{pr description of rank zero c-pairs}, we have $D_c = \SU(n_1) \times \dots \times \SU(n_\ell)$. Let $\delta_{1,i}$ and $\delta_{2,i}$ be the projections of $\delta_1$ and $\delta_2$ to $\SU(n_i)$. The conjugation of the $c$-pair $([s_1, \delta_1]_{\tau_c}, [s_2, \delta_2]_{\tau_c})$ by $[\id, \delta_{1,i}]_{\tau_c}$ gives us $([s_1,\delta_1]_{\tau_c}, [s_2, c_i \delta_2]_{\tau_c})$ and similarly, conjugating by $[\id, \delta_{2,i}]_{\tau_c}$ gives $([s_1, c_i^{-1} \delta_1]_{\tau_c}, [s_2, \delta_2]_{\tau_c})$. By Proposition \ref{pr description of rank zero c-pairs}\eqref{num c is a product of generators of ZZ n}, the $c_i$ generate $Z_{D_c}(D_c)$, so it is obvious that (\ref{eq zeta'}) factors through} 
\beq
\label{eq zeta''}
\ol{S}_c \times \ol{S}_c \lra \Rr(K)_d.
\eeq 
\erm 
One can further prove that (\ref{eq zeta''}) factors through the quotient by the finite group $W_c$, defined in (\ref{eq definition of W_c}).

\btm {\bf (\cite[Corollary 4.2.2]{borel&friedman&morgan}).}
\label{tm homeomorphism RrK}
Let $K$ be a compact connected Lie group. There is a homeomorphism
\[
\quotient{(\ol{S_c} \times \ol{S_c})}{W_c} \stackrel{homeo}{\lra} \Rr(K)_d.
\]
\etm

\brm
\label{rm olzeta}
\nr{Since (\ref{eq definition of Theta}) gives us the isomorphism $\Theta_{\ol{S}_c} :   \U(1) \otimes_\ZZ \Lambda_{\ol{S}_c}\stackrel{\cong}{\lra}\ol{S}_c$ and the action of $W_c$ commutes with $\Theta_{\ol{S}_c}$, we have a natural homeomorphism}  
\[
\quotient{\left( (\U(1) \times \U(1)) \otimes_\ZZ \Lambda_{\ol{S}_c} \right)}{W_c} \stackrel{homeo}{\lra} \Rr(K)_d.
\]
\erm

\subsection{Review on $G$-bundles over an elliptic curve}

Let $G$ be a connected complex reductive Lie group with maximal compact $K$.
The notions of stability, semistability, polystability and S-equivalence for $G$-bundles are well known (see, for example, \cite{ramanathan_stable}). 

Given a unitary representation $\rho : \Gamma_\RR \to K \subset G$, after \cite{atiyah&bott}, we can construct the $G$-bundle $E^\rho$ as follows (see also \cite{ramanathan_stable} for a similar construction). Consider the line bundle $\Oo(x_0)$ associated with the divisor given by the fixed point $x_0$ of $X$ and let $Q'_{x_0} \to X$ be the fixed $\U(1)$-bundle obtained from reduction of structure group of $\Oo(x_0)$. The universal covering $\wt{X} \to X$ is a $\pi_1(X)$-bundle. Consider the fibre product $\wt{X} \times_X Q'_{x_0}$ and denote by $Q_{x_0}$ its lifting to $\Gamma_\RR$. We set $E^\rho$ as the extension of structure group associated to $\rho$ of $Q_{x_0}$, i.e.
\beq
\label{eq definition of E^rho}
E^\rho = \rho_* Q_{x_0}.
\eeq
As shown in \cite{atiyah&bott} (see also \cite{ramanathan_stable}),
\bit
\item the bundles $E^{\rho}$ are polystable,
\item two bundles $E^{\rho_1}$ and $E^{\rho_2}$ are isomorphic if and only if $\rho_1$ and $\rho_2$ are conjugate,
\item every polystable $G$-bundle $E$ is isomorphic to some $E^{\rho}$, and
\item the bundle $E^\rho$ is stable if and only if the representation $\rho$ is irreducible.
\eit 

We can interpret  as follows the results of \cite{borel&friedman&morgan} given in Section \ref{sc borel friedman morgan}. 

\bpr 
\label{pr description of semisimple G such that stable G-bundles exists}
Let $G$ be a connected, complex semisimple Lie group and denote by $\wt{G}$ its universal cover. Let $E^\rho$ be a stable $G$-bundle of degree $d \in \pi_1(G) \subset Z_{\wt{G}}(\wt{G})$. Then
\begin{enumerate}

\item \label{num G is a product of PGLnC} the group $\wt{G}$ is a product of simple factors $G_i$, where each $G_i$ is isomorphic to $\SL(n_i,\CC)$ for some $n_i \geq 2$;

\item \label{num d is a product of generators of ZZ n} $d = (d_1, \dots, d_r)$, where each $d_i$ generates $\pi_1(G_i)\cong \ZZ_{n_i}$;
  
\item \label{num exists a G-bundle as above} conversely, if $G$ is as in (\ref{num G is a product of PGLnC}) and $d$ as in (\ref{num d is a product of generators of ZZ n}), then there is a stable $G$-bundle of degree $d$ and all $G$-bundles of degree $d$ are isomorphic, i.e.
\[
M(G)_d = M^{st}(G)_d = \{ pt \}.
\]
\end{enumerate}
\epr 

\proof
Since, by Remark \ref{rm rho is determined by a b and u}, a representation $\rho$ is determined by a $c$-pair, and the $c$-pair is irreducible if and only if the representation is irreducible, the proof follows from Proposition \ref{pr description of rank zero c-pairs} (\cite[Proposition 4.1.1]{borel&friedman&morgan}) and the existence of a bijective correspondence between irreducible representations and stable $G$-bundles.
\qed

\brm
\nr{Note that Proposition \ref{pr description of semisimple G such that stable G-bundles exists} implies that, for $G$ simple, the only stable bundles occur when $G = \PGL(n,\CC)$ and $d$ generating $\ZZ_n$ ({\it i.e.} $n$ and $d$ coprime).}
\erm

Let $G$ be a complex reductive Lie group, and let $F$ be as defined at the beginning of Section \ref{sc notation and results on Lie groups}. Since $F \subset Z_G(G)$, the extension of structure group given by the multiplication map $\mu : F \times G \to G$ is well defined. Given an $F$-bundle $J$ and a $G$-bundle $E$, we denote by $J \otimes E$ the $G$-bundle $\mu_*(J \times_X E)$. 

\bco
\label{co H^1XF acts trivially on MstG}
Let $E$ be a stable $G$-bundle of topological class $d$ and let $J$ be any element of $H^1(X,F)$. Then
\[
E \cong J \otimes E,
\]
so $J \otimes E$ has the same topological invariant as $E$.
\eco

\proof
This follows from Remark \ref{rm rho conjugated to rho_J rho}.
\qed 

By \cite[Proposition 7.1]{ramanathan_stable}, a $G$-bundle is stable if and only if the induced $(G/Z_0)$-bundle is stable. Let $\ol{Z}$ and $\ol{D}$ be as defined at the beginning of Section \ref{sc notation and results on Lie groups}.

\btm
\label{tm MstG}
Let $G$ be a connected complex reductive Lie group and let $d \in \pi_1(G)$. Then
\[
M^{st}(G)_d = \emptyset,
\]
unless $G/Z_0$ decomposes into $\PGL(n_1,\CC) \times \stackrel{}{\dots} \times \PGL(n_s,\CC)$ and $d\in\pi_1(G_i)$ projects to $(d_1, \dots, d_s) \in \pi_1(\PGL(n_1,\CC)) \times \pi_1(\PGL(n_s,\CC))$ where $\gcd(n_i,d_i) = 1$. In that case, there is a natural isomorphism
\[
M^{st}(G)_d = M(G)_d \cong X \otimes_\ZZ \Lambda_{\ol{Z}}.
\]
\etm  

\proof
The first statement follows from Proposition \ref{pr description of semisimple G such that stable G-bundles exists}.

The extension of structure group associated to $G \to \ol{G} \cong \ol{Z} \times \ol{D}$ (see (\ref{eq olG = olZ times olD})) induces a morphism
\begin{equation} \label{eq M^st G to M^st olG}
M^{st}(G)_d \lra M^{st}(\ol{G})_{(u,c)} \cong M^{st}(\ol{Z})_{u} \times M^{st}(\ol{D})_{c}.
\end{equation}
This morphism is injective by Corollary \ref{co H^1XF acts trivially on MstG}. For any stable $G$-bundle $E$, the morphism
\[
\map{M^{st}(Z)_u}{M^{st}(\ol{Z})_u}{J}{\ol{J} := (J \otimes E)/D \cong J/F}{}
\]
is surjective, as $\ol{J}$ is the extension of structure group of $J$ associated to $Z \to \ol{Z}$. Then, the morphism (\ref{eq M^st G to M^st olG}) is bijective, and, therefore, it is an isomorphism. By Proposition \ref{pr description of semisimple G such that stable G-bundles exists}, $M^{st}(\ol{D})_{c} = \{ pt \}$, so the second statement follows from Theorem \ref{tm MG for G abelian}.
\qed

\brm
\label{rm definition of Ee_G}
\nr{Note that the point $x_0 \in X$ defines an origin in $M^{st}(\ol{Z})_{u}$. For $G$ and $d$ of the form given in Theorem \ref{tm MstG}, we write $E^{x_0}_{G,d}$ for the stable $G$-bundle of degree $d$ associated to this point of $M^{st}(\ol{Z})_{u}$. Let $Z_0$ be the connected component of the centre of $G$ and consider the universal family of $Z_0$-bundles $\Pp_{Z_0}$ pa\-ra\-me\-tri\-zed by $X \otimes_\ZZ \Lambda_{Z_0}$ which is defined in (\ref{eq poincare bundle}). We define the family $(\Ee')_{G,d} = \Pp_{Z_0} \otimes E^{x_0}_{G,d}$ of $G$-Higgs bundles with degree $d$. By Corollary \ref{co H^1XF acts trivially on MstG}, 
 this family 
 descends to a family parametrized 
by the quotient of $X \otimes_\ZZ \Lambda_{Z_0}$ by the image of $H^1(X,F)$. Recalling that $\exp^{-1}(F) = \Lambda_{\ol{Z}} \subset \Lambda_{Z_0}$ one can check that this quotient is isomorphic to $X \otimes_\ZZ \Lambda_{\ol{Z}}$. Then we have a family $\Ee_{G,d} \to X \times (X \otimes_\ZZ \Lambda_{\ol{Z}})$ such that}
\[
\map{X \otimes_\ZZ \Lambda_{\ol{Z}}}{M^{st}(G)_d}{t}{[\Ee_{G,d} |_{X \times \{ t \} } ]_{\cong}.}{\cong}
\]
\erm

\bpr
\label{pr unique Jordan-Holder Levi subgroup for G-bundles}
Every polystable $G$-bundle of topological type $d = (u,c)$ admits a reduction of structure group to $L_c$, giving a stable $L_c$-bundle of topological class $\ell_d = (u,p(c))$.
\epr

\proof
Every polystable $G$-bundle is isomorphic to some $E^\rho$. By Remark \ref{rm rho is determined by a b and u}, $\rho$ is determined by $u$ and a $c$-pair $(a,b) \in K \times K$. By Proposition \ref{pr B xy = B_c and D_c has a unique rk0 c-pair} (\cite[Proposition 4.2.1]{borel&friedman&morgan}), $(a,b)$ is contained (after conjugation) in the maximal compact subgroup of $L_c$ and is irreducible as a $c$-pair in that group. Then $\im \rho \subset L_c$ and $\rho$ is irreducible in $L_c$, so $E^\rho$ reduces to a stable $L_c$-bundle. 
\qed

By Proposition \ref{pr unique Jordan-Holder Levi subgroup for G-bundles}, it makes sense to define the following family parametrizing all polystable $G$-bundles of degree $d$,
\beq
\label{eq definition of Ee_G for arbitrary G}
\Ee_{G,d}:= i_* (\Ee_{L_c,\ell_d}),
\eeq
where $i : L_c \hookrightarrow G$ is the natural inclusion. Note that this family is parametrized by $X \otimes_\ZZ \Lambda_{\ol{S}_c}$, where $S_c$ is the centre of $L_c$. This family induces a morphism to the moduli space
\beq
\label{eq varsigma}
X \otimes_\ZZ \Lambda_{\ol{S}_c} \lra M(G)_d,
\eeq 
which is surjective by Proposition \ref{pr unique Jordan-Holder Levi subgroup for G-bundles}.

\btm 
Let $G$ be a connected complex reductive Lie group and let $d \in \pi_1(G)$. Then 
\[
M(G)_d \cong \quotient{(X \otimes_\ZZ \Lambda_{\ol{S}_c})}{W_c}.
\]
\etm 

\proof
It is clear that (\ref{eq varsigma}) descends to a surjective morphism 
\[
\zeta_{G,d}: \quotient{(X \otimes_\ZZ \Lambda_{\ol{S}_c})}{W_c} \lra M(G)_d.
\]
Injectivity follows from Corollary \ref{co H^1XF acts trivially on MstG} and the fact that the reduction of structure group to $L_c$ is unique up to conjugation. Now $\varsigma_{G,d}$ is an isomorphism by Zariski's Main Theorem. 
\qed

\bco
\label{co varsigma/W_c is injective}
Let $E_1$ and $E_2$ be two polystable $G$-bundles of topological class $d$ pa\-ra\-me\-tri\-zed by $\Ee_{G,d}$ at the points $t_1$ and $t_2 \in X \otimes_\ZZ \Lambda_{\ol{S}_c}$. Then $E_1$ and $E_2$ are isomorphic $G$-bundles if and only if there exists $\omega \in W_c$ such that $t_2 = \omega \cdot t_1$.
\eco

\section{$G$-Higgs bundles over an elliptic curve}
\label{sc chapter G-Higgs bundles}

Let $G$ be a connected complex reductive Lie group. Recall that a $G$-Higgs bundle over an elliptic curve $X$ is a pair $(E,\Phi)$, where $E$ is an algebraic $G$-bundle over $X$ and $\Phi\in H^0(X,E(\gG))$. 
We say that $(E,\Phi)$ is {\it stable} (resp. {\it semistable}) if, for every proper parabolic subgroup $P$ with Lie algebra $\pP$, any non-trivial antidominant character $\chi: P \to \CC^*$, and any reduction of structure group $\sigma$ to the parabolic subgroup $P$ giving the $P$-bundle $E_{\sigma}$ such that $\Phi \in H^0(X,E_\sigma(\pP))$, we have
\[
\deg \chi_* E_{\sigma} > 0 \quad (\nr{resp.} \qua \geq 0).
\]

Let $(E_1,\Phi_1)$ and $(E_2, \Phi_2)$ be two semistable $G$-Higgs bundles and suppose that there exists a family $\Hh$ parametrized by $\CC$ such that $\Hh|_{X \times \{ \lambda \}} \cong (E_1,\Phi_1)$ if $\lambda \neq 0$ and $\Hh|_{X \times \{ 0 \}} \cong (E_2,\Phi_2)$. We say that these two $G$-Higgs bundles are {\it S-equivalent} and we call the induced equivalence relation {\it S-equivalence}, writing $(E_1,\Phi_1) \sim_S (E_2,\Phi_2)$. Two families of semistable $G$-Higgs bundles parametrized by $Y$ are {\it S-equivalent}, $\Hh_1 \sim_S \Hh_2$, if for every point $y \in Y$, one has $\Hh_1|_{X \times \{ y \}} \sim_S \Hh_2 |_{X \times \{ y \}}$. 

We denote by $\Mm(G)_d$ the moduli space of S-equivalence classes of semistable $G$-Higgs bundles of degree $d$ and by $\Mm^{st}(G)_d$ the corresponding moduli space for stable $G$-Higgs bundles.

The $G$-Higgs bundle $E$ is {\it polystable} if it is semistable and, when there exists a parabolic subgroup $P \subsetneq G$, a strictly antidominant character $\chi : P \to \CC^*$ and a reduction of structure group $\sigma$ giving the $P$ bundle $E_\sigma$ such that
\[
\Phi \in H^0(X,E_{\sigma}(\pP))
\]
and
\[
\deg \chi_* E_{\sigma}  = 0,
\]
there exists a reduction $\varsigma$ of the structure group of $E_\sigma$ to the Levi subgroup $L \subset P$ such that $\Phi \in H^0(X,E_{\varsigma}(\lL))$, where $E_\varsigma$ denotes the principal $L$-bundle obtained from the reduction of structure group $\varsigma$ and $\lL$ is the Lie algebra of $L$. There is a unique (up to isomorphism) polystable $G$-Higgs bundle in each S-equivalence class. Let us recall that every polystable $G$-Higgs bundle has a reduction of structure group to some Levi subgroup $L \subset G$ giving a stable $L$-Higgs bundle. Such a reduction is called a {\it Jordan--Hölder reduction} and is unique in a certain sense (see, for example, \cite{oscar&ignasi&gothen}).

The triviality of the canonical bundle $\Omega_X^1$ in the case of an elliptic curve leads us to the following well known results.

\bpr
\label{pr E Phi semistable iff E semistable for G-Higgs bundles}
Let $(E,\Phi)$ be a semistable $G$-Higgs bundle. Then $E$ is a semistable $G$-bundle.
\epr 

\proof
If $E$ is unstable, then $E$ reduces to the Harder--Narasimhan parabolic subgroup $P$, giving $E_\sigma$, and there exists a character $\chi : P \to \CC^*$ such that $\deg \chi^* E_\sigma < 0$. Moreover $H^0(X,E(\lie{g})) = H^0(X,E_\sigma(\lie{p}))$. So $\Phi \in H^0(X,E_\sigma(\lie{p}))$ and hence the Higgs bundle $(E,\Phi)$ is unstable.
\qed

We have the following consequence.

\bco
\label{co MmG onto MG}
The moduli space of $G$-Higgs bundles projects onto the moduli space of $G$-bundles
\[
\map{\Mm(G)_{d}}{M(G)_{d}}{\left[(E,\Phi)\right]_S}{\left[E\right]_S.}{}
\]
\eco

\bpr
\label{pr E Phi stable iff E stable for G-Higgs bundles}
Let $(E,\Phi)$ be a stable $G$-Higgs bundle. Then $E$ is stable.
\epr

\proof
We first note that $\Phi \in H^0(X,E(\gG))$ is contained in $\aut(E,\Phi)$.

If $(E,\Phi)$ is stable, then, by \cite[Proposition 2.14]{oscar&ignasi&gothen}, $\aut(E,\Phi) \subset H^0(X,E(\zZ_{\gG}(\gG)))$ and it follows easily that $(E,0)$ is stable too. 
\qed

\bco
\label{co E Phi polystable implies E polystable}
Let $(E,\Phi)$ be a polystable $G$-Higgs bundle. Then $E$ is a polystable $G$-bundle.
\eco

\proof
The polystable $G$-Higgs bundle $(E,\Phi)$ reduces to the Jordan--Hölder Levi subgroup $L$ giving the stable $L$-Higgs bundle $(E_L,\Phi_L)$. By Proposition \ref{pr E Phi stable iff E stable for G-Higgs bundles}, $E_L$ is a stable $L$-bundle and therefore $E$ is a polystable $G$-bundle.
\qed 

With the results above we are able to describe stable and polystable $G$-Higgs bundles. Recall the bundle $E^\rho$ defined in (\ref{eq definition of E^rho}). 

\bpr
\label{pr description stable G-Higgs bundles}
A stable $G$-Higgs bundle $(E,\Phi)$ is isomorphic to $(E^{\rho},z \otimes 1_{\Oo})$ where $\rho : \Gamma_\RR \to K$ is some representation such that $\zZ_\gG(\rho) = \zZ_{\gG}(\gG)$, $1_{\Oo}$ is the constant section of the trivial bundle $\Oo$ equal to $1$ and $z \in \zZ_{\gG}(\gG)$. 
\epr

\proof
By Proposition \ref{pr E Phi stable iff E stable for G-Higgs bundles}, $E$ is stable and therefore polystable. Then $E \cong E^\rho$ for some $\rho$. By \cite[Proposition 3.2]{ramanathan_stable}, we have $H^0(X,E(\gG)) = \zZ_{\gG}(\gG)$, so $\Phi = z \otimes 1_{\Oo}$ for some $z \in \zZ_{\gG}(\gG)$. Note that $\zZ_{\gG}(\rho) \subseteq H^0(X,E^{\rho}(\gG))$, and then $\zZ_{\gG}(\rho)$ is contained in $\zZ_{\gG}(\gG)$ so they are equal.
\qed

We recall the isomorphism (\ref{eq definition of dTheta}) and note that $T^*X \cong X \times \CC$. With all this in mind, we provide a result for $G$-Higgs bundles analogous to Theorem \ref{tm MG for G abelian}. 
 
\btm \label{tm Mmst fo G abelian}
Let $S^\CC$ be an abelian, connected complex Lie group. Then, the moduli space of topologically trivial $S^\CC$-Higgs bundles over the elliptic curve $X$ is
\[
\Mm(S^\CC)_0 \cong T^*X \otimes_\ZZ \Lambda_{S}.
\]
\etm 

\proof
The description follows from the construction of a family of $S^\CC$-Higgs bundles using $\Pp_S$ defined in (\ref{eq poincare bundle}) and $d \Theta_{S}$ from  (\ref{eq definition of dTheta}).
\qed

Recall the definition of $L_c$ given in (\ref{eq definition of L_c}). 

\bpr
\label{pr unique Jordan-Holder Levi subgroup for G-Higgs bundles}
Every polystable $G$-Higgs bundle of topological type $d = (u,c)$ admits a reduction of structure group to $L_c$ giving a stable $L_c$-Higgs bundle of topological class $\ell_d = (u,p(c))$.
\epr

\proof
Take a polystable $G$-Higgs bundle $(E,\Phi)$ of type $d=(u,c)$, and suppose that $L$ is a Jordan--Hölder Levi subgroup of $(E,\Phi)$. Since $(E,\Phi)$ reduces to $L$ giving a stable $L$-Higgs bundle, it follows from Proposition \ref{pr description stable G-Higgs bundles} that there exists $(\rho,z)$ such that $(E,\Phi) \cong (E^\rho,z \otimes 1_\Oo)$. Here $z \in \zZ_{\gG}(\rho)$, which is a reductive Lie algebra since 
\[
Z_{G}(\rho) = Z_{G}(a,b) = Z_K(a,b)^\CC,
\]
and $Z_K(a,b)$ is a compact subgroup. Then we can conjugate $z \in \zZ_\gG(\rho)$ to the Cartan subalgebra $\hH \subset \lL_c$. As a consequence of the above and Proposition \ref{pr unique Jordan-Holder Levi subgroup for G-bundles}, $(E^\rho,z\otimes1_\Oo)$ reduces to a stable $L_c$-Higgs bundle and so does $(E,\Phi)$. 
\qed

Recall that $\hH^{\omega_c}$ is the centre of $\lL_c$. Propositions \ref{pr description stable G-Higgs bundles} and \ref{pr unique Jordan-Holder Levi subgroup for G-Higgs bundles} imply the following.

\bco
\label{co description polystable G-Higgs bundles}
Let $K_{L_c}$ be a maximal compact subgroup of $L_c$. A polystable $G$-Higgs bundle $(E,\Phi)$ of type $d \in \pi_1(G)$ is isomorphic to $(E^{\rho},z \otimes 1_{\Oo})$ where $\rho : \Gamma_\RR \to K_{L_c} \subset L_c$ is some representation, $1_{\Oo}$ is the constant section of the trivial bundle $\Oo$ equal to $1$ and $z \in \hH^{\omega_c}$.
\eco

Recall that $E^{\rho_1} \cong E^{\rho_2}$ if and only if $\rho_1$ and $\rho_2$ are conjugate. This fact, together with Corollary \ref{co description polystable G-Higgs bundles}, implies the following.

\bco
\label{co automorphisms groupn of polystable G-Higgs bundles}
In the notation of Corollary \ref{co description polystable G-Higgs bundles}, two pairs $(\rho, z)$ and $(\rho', z')$ determine isomorphic polystable $G$-Higgs bundles if and only if there exists an element $k \in K$ such that $(\rho', z') = (k \rho k^{-1}, \ad_k(z))$.

The automorphism group of the polystable $G$-Higgs bundle $(E^\rho, z \otimes 1_\Oo)$ is $Z_G(\rho,z)$ and its Lie algebra is $\zZ_{\gG}(\rho,z)$.
\eco

Recall the family of polystable $G$-bundles $\Ee_{G,d} \to X \times (X \otimes_\ZZ \Lambda_{\ol{S}_c})$ defined in Remark \ref{rm definition of Ee_G} and in (\ref{eq definition of Ee_G for arbitrary G}). Recalling the isomorphism 
\[
d \Theta_{\ol{S}_c} : \CC \otimes_\ZZ \Lambda_{\ol{S}_c} \to \sS_c = \hH^{\omega_c} 
\]
defined in (\ref{eq definition of dTheta}), as well as the discussion immediately before Theorem \ref{tm Mmst fo G abelian}, we define a family of $G$-Higgs bundles $\Hh_{G,d}$ parametrized by $T^*X \otimes_\ZZ \Lambda_{\ol{S}_c}$, setting, for each point $(t,s) \in T^*X \otimes_\ZZ \Lambda_{\ol{S}_c}$, 
\[
\Hh_{G,d} |_{X \times \{ (t,s) \}} = \left( \Ee_{G,d} |_{X \times \{ t \}} \qua , \qua d\Theta_{\ol{S}_c}(s) \otimes 1_\Oo  \right), 
\]
where $1_\Oo$ is the section of the trivial bundle $\Oo$ equal to $1$.

\brm
\label{rm eta is surjective}
\nr{By Corollary \ref{co description polystable G-Higgs bundles}, every polystable $G$-Higgs bundle of degree $d$ is pa\-ra\-me\-tri\-zed by $\Hh_{G,d}$.}
\erm

\brm
\label{rm abelianization}
\nr{The family $\Hh_{G,d}$ can be constructed starting from $\Hh_{S_c,0} \otimes (E^{x_0}_{L_c,\ell_d},0)$, quotienting by $H^1(X,F)$ as described in Corollary \ref{co H^1XF acts trivially on MstG} and taking the extension of structure group associated to $L_c \hookrightarrow G$. This shows that all polystable $G$-Higgs bundles are described by Higgs bundles for the abelian group $S_c$.}
\erm

\btm
\label{tm MmMstG}
Let $G$ be a connected complex reductive Lie group and let $d \in \pi_1(G)$. Then
\[
\Mm^{st}(G)_d = \emptyset,
\]
unless $G/Z_0$ decomposes into $\PGL(n_1,\CC) \times \stackrel{}{\dots} \times \PGL(n_s,\CC)$ and $d\in\pi_1(G)$ projects to $(d_1, \dots, d_s) \in \pi_1(\PGL(n_1,\CC)) \times \pi_1(\PGL(n_s,\CC))$ where $\gcd(n_i,d_i) = 1$. In that case, there is a natural isomorphism
\[
\Mm^{st}(G)_d = \Mm(G)_d \cong T^*X \otimes_\ZZ \Lambda_{\ol{Z}}.
\]
\etm 

\proof
The first statement is a consequence of Propositions \ref{pr E Phi stable iff E stable for G-Higgs bundles} and \ref{pr description stable G-Higgs bundles} and Theorem \ref{tm MstG}.

As in Theorem \ref{tm MstG}, the extension of structure group associated to $G \to \ol{G} \cong \ol{Z} \times \ol{D}$ induces a morphism
\begin{equation} \label{eq Mm^st G to Mm^st olG}
\Mm^{st}(G)_d \lra \Mm^{st}(\ol{G})_{(u,c)} \cong \Mm^{st}(\ol{Z})_{u} \times \Mm^{st}(\ol{D})_{c},
\end{equation}
which, as in the case of $G$-bundles, can be proved to be bijective. By Proposition \ref{pr description stable G-Higgs bundles}, Corollary \ref{co automorphisms groupn of polystable G-Higgs bundles} and Theorem \ref{tm MstG}, $\Mm^{st}(\ol{D})_{c} = \{ pt \}$. Noting also that $\Mm^{st}(\ol{Z})_{u}$ is smooth as $\ol{Z}$ is abelian, we have that $\Mm^{st}(G)_d \cong \Mm^{st}(\ol{Z})_{u}$ and the second statement follows from Theorem \ref{tm Mmst fo G abelian}.
\qed

Recall $W_c$ defined in (\ref{eq definition of W_c}). Note that $W_c$ acts on $\ol{S}_c$ and therefore it acts on $T^*X \otimes_\ZZ \Lambda_{\ol{S}_c}$.

\bpr
\label{pr eta/W_c is injective}
Let $(E_1,\Phi_1)$ and $(E_2,\Phi_2)$ be two polystable $G$-Higgs bundles of topological class $d$ parametrized by $\Hh_{G,d}$ at the points $(t_1,s_1)$ and $(t_2,s_2) \in T^*X \otimes_\ZZ \Lambda_{\ol{S}_c}$. Then $(E_1,\Phi_1)$ and $(E_2,\Phi_2)$ are isomorphic $G$-Higgs bundles if and only if there exists $\omega' \in W_c$ such that $(t_2,s_2) = \omega' \cdot (t_1,s_1)$.
\epr  

\proof
It is clear that, if $(t_2,s_2) = \omega' \cdot (t_1,s_1)$, then $(E_1,\Phi_1) \cong (E_2,\Phi_2)$. Suppose conversely that $(E_1,\Phi_1) \cong (E_2,\Phi_2)$ and that $(E_1,\Phi_1)$ and $(E_2,\Phi_2)$ are associated to $(\rho_1,z_1)$ and $(\rho_2,z_2)$ in the sense of Corollary \ref{co description polystable G-Higgs bundles}. Then, by Corollary \ref{co automorphisms groupn of polystable G-Higgs bundles}, there exists $k \in K$ such that $(\rho_2,z_2)$ is equal to $(k\rho_1 k^{-1}, \ad_k z)$.

By Corollary \ref{co varsigma/W_c is injective}, there exists $\omega \in W_c = N_K(S_c)/Z_K(S_c)$ such that $t_2 = \omega \cdot t_1$. Then there exists $n \in N_K(L_c)=N_K(S_c)$ projecting to $\omega$ and such that $\rho_2 = n \rho_1 n^{-1}$. Let us set $z' = \ad_{n^{-1}} (z_2)$ in $\sS_c = \hH^{\omega_c}$ and note that
\[
(\rho_2,z_2) =  \left( n \rho_1 n^{-1}, \ad_n(z') \right).
\]
Then $(\rho_1, z') = \left( (n^{-1}k) \rho_1 (n^{-1}k)^{-1}, \ad_{n^{-1}k} z_1 \right)$, so $n^{-1}k$ belongs to $Z_K(\rho_1)$ and conjugates $z_1$ to $z'$, both elements of $\sS_c = \hH^{\omega_c}$. 

Let $T$ be the maximal torus of $Z_K(\rho_1,z')$ such that its complexification is $S_c$. Note that $T' = n^{-1}k T (n^{-1}k)^{-1}$ is another maximal torus of $Z_K(\rho_1,z')$. Since $Z_K(\rho_1,z')$ is compact there exists an element $h'$ that conjugates $T$ to $T'$. Then, there exists $h = n^{-1}kh' \in Z_K(\rho_1) \cap N_K(S_c)$ with $z' = \ad_{h}(z_1)$. Setting $n' = n h = kh'$ we obtain an element of $N_K(S_c)$ such that 
\[
(\rho_2, z_2) = \left( n' \rho_1 (n')^{-1}, \ad_{n'} (z_1) \right).
\]

Finally, let $\omega' \in W_c$ be given by the projection of $n'$. It is clear that it sends $(t_1,s_1)$ to $(t_2,s_2)$.
\qed

\btm
\label{tm MmG is the normalization of MmMG}
There exists a bijective morphism 
\beq 
\label{eq eta_G d}
\quotient{(T^*X \otimes_\ZZ \Lambda_{\ol{S}_c})}{W_c} \stackrel{1:1}{\lra} \Mm(G)_d.
\eeq
Hence the normalization $\ol{\Mm(G)}_d$ of $\Mm(G)_d$ is isomorphic to $\quotient{(T^*X \otimes_\ZZ \Lambda_{\ol{S}_c})}{W_c}$.
\etm

\bpf
By moduli theory, the family $\Hh_{G,d} \to X \times (T^*X \otimes_\ZZ \Lambda_{\ol{S}_c})$ induces a morphism
\[
\map{T^*X \otimes_\ZZ \Lambda_{\ol{S}_c}}{\Mm(G)_d}{(t,s)}{[\Hh_{G,d} |_{X \times \{ (t,s) \}}]_S.}{}
\]
As we have seen in Remark \ref{rm eta is surjective}, this morphism is surjective. It descends to a surjective morphism (\ref{eq eta_G d}). By Proposition \ref{pr eta/W_c is injective}, (\ref{eq eta_G d}) is also injective.

The quasiprojective variety $\quotient{(T^*X \otimes_\ZZ \Lambda_{\ol{S}_c})}{W_c}$ is normal since it is the quotient of a smooth (and therefore normal) variety by a finite (and therefore reductive) group. Zariski's Main Theorem and (\ref{eq eta_G d}) give us the description of the normalization of $\Mm(G)_d$.
\epf

\brm
\label{rm Normality of MmG}
\nr{This is proved in \cite{thaddeus} for the trivial degree case. For $G = \GL(n,\CC)$ or $\SL(n,\CC)$ and $d = 0$, (\ref{eq eta_G d}) is indeed an isomorphism since the target is normal (see the discussion at the end of Section \ref{sc introduction}).}
\erm 

The irreducibility of the quotient $\quotient{(T^*X \otimes_\ZZ \Lambda_{\ol{S}_c})}{W_c}$ implies the following.

\bco
The moduli space of $G$-Higgs bundles $\Mm(G)_d$ is irreducible.
\eco

A $G$-Higgs bundle is {\it infinitesimally regular} if the dimension of $\aut(E,\Phi)$ is the minimal possible one.

\bpr
\label{pr smooth points of MmMG_d} 
The Zariski open subset of points represented by polystable $G$-Higgs bundles which are infinitesimally regular lies in the smooth locus of $\Mm(G)_d$. 
\epr

\proof
Consider the infinitesimal deformation space $T$ of $(E,\Phi)$. By \cite{biswas&ramanan} one has the exact sequence
\[
H^0(X,E(\gG))\stackrel{e_0(\Phi)}{\lra} H^0(X,E(\gG)) \lra T\lra H^1(X,E(\gG)) \stackrel{e_1(\Phi)}{\lra} H^1(X,E(\gG)),
\]
where $e_i(\Phi)(\psi) = [\psi, \Phi]$ and $e_1(\Phi)$ is the Serre dual of $e_0(\Phi)$ (recall that the canonical bundle is trivial in our case). Hence $\codim(\im e_0(\Phi)) = \dim(\ker e_1(\Phi))$, so $\dim(T)=2 \dim(\ker e_1(\Phi))$.

Suppose that $(E,\Phi) \cong (E^\rho, z \otimes 1_\Oo)$. Recall that $dx$ is a generator of $H^1(X,\Oo)$, so $H^1(X,E(\gG)) = \{ z' \otimes dx : z' \in \zZ_{\gG}(\rho) \}$. We observe that the kernel of $e_1(\Phi)$ corresponds to $\zZ_{\gG}(\rho,z)$ and therefore
\[
\dim(T)=2 \dim(\zZ_{\gG}(\rho,z)) = 2 \dim (\aut(E,\Phi)),
\] 
where the last step in the equality follows from Corollary \ref{co automorphisms groupn of polystable G-Higgs bundles}.

Suppose that $\rho$ is associated to the $c$-pair $(a,b)$ with (up to conjugation) $a \in H$. Recall that Proposition \ref{pr B xy = B_c and D_c has a unique rk0 c-pair} implies that  $\hH^{\omega_c}$ is a Cartan subalgebra of $\zZ_{\gG}(\rho)$ and therefore a Cartan subalgebra of $\zZ_{\gG}(\rho,z)$ since $z \in \hH^{\omega_c}$. Then, for every polystable $G$-Higgs bundle $(E,\Phi)$,
\[
\dim(\Mm(G)_d) = 2 \dim(\hH^{\omega_c}) \leq 2 \dim(\zZ_{\gG}(\rho,z)) = 2 \dim(\aut(E,\Phi)).
\]
Recalling \cite[Corollary 5.18]{friedman&morgan}, we observe that, if $a$ is a regular element of $H$, then $\zZ_{\gG}(\rho,z) = \hH^{\omega_c}$, so $\dim(T) = \dim(\Mm(G)_d)$ is achieved in a Zariski open subset and the statement follows.
\qed

We define the projection
\[
\morph{p_{G,d}}{\quotient{(T^*X \otimes_\ZZ \Lambda_{\ol{S}_c})}{W_c}}{\quotient{(X \otimes_\ZZ \Lambda_{\ol{S}_c})}{W_c}}{[(t,s)]_{W_c}}{[t]_{W_c}.}{}
\]

Recalling the projection of Corollary \ref{co MmG onto MG}, we have the commutative diagram
\[
\xymatrix{
\quotient{(T^*X \otimes_\ZZ \Lambda_{\ol{S}_c})}{W_c} \ar[d]_{1:1} \ar[rr]^{p_{G,d}} & & \quotient{(X \otimes_\ZZ \Lambda_{\ol{S}_c})}{W_c} \ar[d]^{\cong} 
\\
\Mm(G)_d \ar[rr] & &  M(G)_d.
}
\]

\brm
\label{rm orbifold structure}
\nr{We can give an interpretation of the projection $p_{G,d}$
 in terms of a certain orbifold bundle. Given an orbifold defined as a global quotient $Z/\Gamma$, one can define its cotangent orbifold bundle as the orbifold given by $T^*Z/\Gamma$, where the action of $\Gamma$ on $T^*Z$ is the action induced by the action of $\Gamma$ on $Z$. Denote by $\wt{M}(G)_d$ and $\wt{\Mm}(G)_d$ the orbifolds given respectively by the quotients of $(X \otimes_\ZZ \Lambda_{\ol{S}_c})$ and $(T^*X \otimes_\ZZ \Lambda_{\ol{S}_c})$ by the finite group $W_c$. Since $T^*(X \otimes_\ZZ \Lambda_{\ol{S}_c})$ is $(T^*X \otimes_\ZZ \Lambda_{\ol{S}_c})$, we have that $\wt{\Mm}(G)_d$ is the cotangent orbifold bundle of $\wt{M}(G)_d$, i.e.}
\[
\wt{\Mm}(G)_d \cong \Tt^* \wt{M}(G)_d.
\]
\erm

\section{The Hitchin fibration}
\label{sc the hitchin fibration}

We describe the Hitchin map in the spirit of \cite{donagi&pantev}. Consider the adjoint action of the group $G$ on the Lie algebra $\lie{g}$ and take the quotient map 
\[
q : \lie{g} \lra \git{\gG}{G}.
\]
Let $E$ be any holomorphic $G$-bundle. Since the adjoint action of $G$ on $\lie{g} \Slash G$ is obviously trivial, we note that the fibre bundle induced by $E$ is trivial
\[
E(\lie{g} \Slash G) = \Oo \otimes (\gG \Slash G).
\]
The projection $q$ induces a surjective morphism of fibre bundles
\[
q_E : E(\gG) \lra E(\gG \Slash G),
\]
and $q_E$ induces a morphism on the set of holomorphic global sections
\[
\morph{(q_E)_*}{H^0(X,E(\gG))}{H^0(X,\Oo \otimes (\gG \Slash G))}{\Phi}{\Phi \Slash G}{}.
\]

If $(E_1, \Phi_1)$ and $(E_2,\Phi_2)$ are two S-equivalent semistable $G$-Higgs bundles, one can check that $(q_{E_1})_* \Phi_1 = (q_{E_2})_* \Phi_2$. Hence we can define the Hitchin map
\beq
\label{eq definition of b_G}
\morph{b_G}{\Mm(G)}{H^0(X,\Oo \otimes (\gG \Slash G))}{[(E,\Phi)]_S}{(q_E)_* \Phi.}{}
\eeq
When the base variety is a Riemann surface of genus greater than or equal to $2$, the restriction of $b_G$ to every component $\Mm(G)_d$ is surjective. This is not the case for genus $g = 1$ and, to preserve the fact that the Hitchin map is a fibration, we set
\[
B(G,d) := b_G(\Mm(G)_d),
\]
and we denote by $b_{G,d}$ the restriction of (\ref{eq definition of b_G}) to $\Mm(G)_d$.

If $H$ is a Cartan subgroup with Cartan subalgebra $\hH$ and Weyl group $W$, Chevalley's Theorem says that 
\[
\gG \Slash G \cong \hH / W. 
\]
So $H^0(X, \Oo \otimes (\gG \Slash G)) \cong H^0(X, \Oo \otimes \hH / W)$ and, since $X$ is a compact holomorphic variety, we have $H^0(X,\Oo \otimes \hH/W) \cong \hH/W \cong \CC \otimes_\ZZ \Lambda_H/W$. There is a natural isomorphism
\[
\beta_{G,0}: \quotient{(\CC \otimes_\ZZ \Lambda_H)}{W} \stackrel{\cong}{\lra} B(G,0).
\]

Now we take $d \in \pi_1(G)$ non-trivial associated to $(u,c) \in \pi_1(\ol{Z}) \times \pi_1(\ol{D})$. By Corollary \ref{co description polystable G-Higgs bundles} we see that every polystable $G$-Higgs bundle of topological class $d$ is isomorphic to $(E^{\rho},z \otimes 1_\Oo)$ where $z \in \hH^{\omega_c}$. We can check that the quotient map $q$ induces a bijective morphism
\[
\beta_{G,d} : \quotient{(\CC \otimes_\ZZ \Lambda_{\ol{S}_c})}{W_c} \stackrel{1:1}{\lra} B(G,d).
\]

Let $\Bb(\Lambda_{\ol{S}_c}) = \{ \gamma_1, \dots, \gamma_\ell \}$ be a basis of $\Lambda_{\ol{S}_c}$. Recalling that $T^*X \cong X \times \CC$, we see that the projection $\pi : T^*X \to \CC$ induces
\beq
\label{eq definition of pi_G c}
\morph{\pi_{G,c}}{\quotient{(T^*X \otimes_{\ZZ} \Lambda_{\ol{S}_c})}{W_c}}{\quotient{(\CC \otimes_{\ZZ} \Lambda_{\ol{S}_c})}{W_c}}{[(t,s)]_{W_c}}{[s]_{W_c}.}{}
\eeq
We use this morphism to better understand the Hitchin map.

\bpr
\label{pr description of the Hitchin fibration}
Recall the bijective morphism (\ref{eq eta_G d}). The following diagram is commutative:
\[
\xymatrix{
\quotient{(T^*X \otimes_{\ZZ} \Lambda_{\ol{S}_c})}{W_c}  \ar[rr]^{1:1} \ar[d]_{\pi_{G,c}} & &  \Mm(G)_d \ar[d]^{b_{G,d}}
\\
\quotient{(\CC \otimes_{\ZZ} \Lambda_{\ol{S}_c})}{W_c}  \ar[rr]^{\beta_{G,d}}_{1:1} & & B(G,d).
}
\]
The normalization of the Hitchin fibre corresponding to $s \in \CC \otimes_\ZZ \Lambda_{S_c}$ is isomorphic to
\beq
\label{eq description of pi^-1 of s}
\pi^{-1}_{G,c}([s]_{W_c}) \cong \quotient{(X \otimes_\ZZ \Lambda_{\ol{S}_c})}{Z_{W_c}(s)}.
\eeq 
\epr

\proof
Take $(t,s) \in (T^*X \otimes_{\ZZ} \Lambda_{\ol{S}_c})$, and consider
\[
b_G(\Hh_{G,d}|_{X \times \{(t,s)\}}) = [s]_G.
\]
Clearly, this equality is $W_c$-invariant. On the other hand, note that
\[
\beta_{G,d} \circ \pi_{G,c}\left([(t,s)]_{W_c}\right) = \beta_{G,d}([s]_{W_c}) = [s]_G 
\]
and the first statement follows.

Next, consider the following projection
\[
\wt{\pi}_{G,c} : T^*X \otimes_\ZZ \Lambda_{\ol{S}_c} \lra \CC \otimes_\ZZ \Lambda_{\ol{S}_c}.
\]
We observe that 
\[
\pi^{-1}_{G,c} ([s]_{W_c}) \cong \quotient{\left( \bigcup_{\omega \in W_c} \wt{\pi}_{G,c}^{-1}(\omega \cdot s)\right)}{W_c}.
\]
Since, for $\omega \cdot s \neq \omega' \cdot s$ the sets $\wt{\pi}_{G,c}^{-1}(\omega \cdot s)$ and $\wt{\pi}_{G,c}^{-1}(\omega' \cdot s)$ are disjoint, it follows that 
\[
\pi^{-1}_{G,c}([s]_{W_c}) \cong \quotient{\wt{\pi}_{G,c}^{-1}(s)}{Z_{W_c}(s)}
\]
and therefore we obtain the isomorphism (\ref{eq description of pi^-1 of s}). Finally we observe that the bijection (\ref{eq eta_G d}) sends $\pi^{-1}_{G,c}([s]_{W_c})$ to the Hitchin fibre corresponding with the Higgs field $\Phi = z \otimes 1_\Oo$. Hence, by Zariski's Main Theorem, it describes an isomorphism with the normalization of this subset.
\qed

We denote by $U_{G,c}$ the subset of $\CC \otimes_\ZZ \Lambda_{\ol{S}_c}/W_c$ given by the points $[s]_{W_c}$ such that there exists a non-trivial $\omega \in W_c$ with $s = \omega \cdot s$. Since the only element of $W_c$ that acts trivially on $\CC \otimes_\ZZ \Lambda_{\ol{S}_c}$ is the identity, $U_{G,c}$ is a finite union of closed subsets of codimension at least equal to $1$. By construction, for any $s \notin U_{G,c}$ we have $Z_{W_c}(s) = \{ \id \}$. 

The {\it generic Hitchin fibre} is the fibre over any element of the complement of $U_{G,c}$.

\bco
\label{co generic Hitchin fibre}
The normalization of the generic Hitchin fibre is isomorphic to the abelian variety $X \otimes_{\ZZ} \Lambda_{\ol{S}_c}$.
\eco

\section{The moduli space of representations $\Rr(G)_d$}
\label{sc representations}

From the non-abelian Hodge correspondence on a compact Riemann surface \cite{hitchin-self_duality_equations, simpson2, donaldson, corlette}, it follows that a polystable $G$-Higgs bundle is associated to a reductive representation $\rho : \Gamma_\RR \to G$ and two representations are conjugate if and only if they are associated to isomorphic polystable $G$-Higgs bundles. Furthermore, irreducible representations correspond to stable $G$-Higgs bundles.

Using this correspondence and Remark \ref{rm rho is determined by a b and u}, we can use the results on $G$-Higgs bundles obtained in Section \ref{sc chapter G-Higgs bundles}, to generalize the description given in Section \ref{sc representations and c-pairs} of $c$-pairs on compact groups, to complex reductive Lie groups.

\bpr  
\label{co description of rank zero c-pairs for G}
Let $G$ be a simply connected complex semisimple Lie group. Let $C = p^{-1}(F) = \pi_1(\ol{D})$ as defined at the begining of Section \ref{sc notation and results on Lie groups} and set $c \in C$. Let $(a,b)$ be an irreducible $c$-pair in $G$. Then
\begin{enumerate}

\item \label{num G is a product of SLnC} the group $G$ is a product of simple factors $G_i$, where each $G_i$ is isomorphic to $\SL(n_i,\CC)$ for some $n_i \geq 2$;

\item \label{num c is a product of generators of ZZ n for G} $c = (c_1, \dots, c_r)$, where each $c_i$ generates the centre of $G_i$;
  
\item \label{num a c-pair as above has rank zero for G} conversely, if $G$ is as in (\ref{num G is a product of SLnC}) and $c$ as in (\ref{num c is a product of generators of ZZ n for G}), then there is an irreducible $c$-pair in $G$ and all $c$-pairs in $G$ are conjugate.

\end{enumerate}
\epr

\proof
This follows from Theorem \ref{tm MmMstG} and the fact that the universal cover of $\PGL(n,\CC)$ is $\SL(n,\CC)$.
\qed

\bpr 
\label{co B xy = B_c and D_c has a unique rk0 c-pair for G}
Let $G$ be a connected complex reductive Lie group. Let $(a,b)$ be a reductive $c$-pair; then $(a,b)$ is contained in $L_c$ after conjugation and, as a $c$-pair in $L_c$, is irreducible. 
\epr 

\proof
This follows from Proposition \ref{pr unique Jordan-Holder Levi subgroup for G-Higgs bundles}.
\qed 

Recall the notation introduced in Section \ref{sc notation and results on Lie groups}.

\btm
\label{tm homeomorphism RrG}
Let $G$ be a connected complex reductive Lie group and let $d \in \pi_1(G)$, corresponding under the injection (\ref{eq definition of q^pi}) to $(u,c) \in \pi_1(\ol{Z}) \times \pi_1(\ol{D})$. Then there is a bijective morphism
\[
\zeta_{G,d}: \quotient{(\ol{S_c} \times \ol{S_c})}{W_c} \stackrel{1:1}{\lra} \Rr(G)_d.
\]
\etm

\proof
Take a representative $(\delta_1, \delta_2)$ of the unique conjugation class of $c$-pairs in $D_c$. Recall that $C(G)^+_c$ denotes the space of reductive $c$-pairs in $G$ and consider the following morphisms
\[
\mapp{(S_c \times S_c)}{C(G)^+_c}{\Rr(G)_d}{(s_1,s_2)}{([s_1, \delta_1]_{\tau_c}, [s_2, \delta_2]_{\tau_c})}{\left[ ([s_1, \delta_1]_{\tau_c}, [s_2, \delta_2]_{\tau_c}) \right]_G.}{}
\]
By an argument analogous to that of Remark \ref{rm rho conjugated to rho_J rho}, the composition morphism factors through 
\[
\ol{S}_c \times \ol{S}_c \lra \Rr(G)_d.
\] 
By (\ref{eq reductive c-pairs are enougth}) and Proposition \ref{co B xy = B_c and D_c has a unique rk0 c-pair for G}, it is clear that this morphism is surjective. The group $W_c$ acts on $\ol{S}_c \times \ol{S}_c$ via conjugation by $N_G(S_c)$. Since the points of $\Rr(G)_d$ are the conjugation classes of $c$-pairs, the morphism factors through this quotient, giving the morphism $\zeta_{G,d}$ of the statement. We only need to prove that it is injective. 

Take two reductive $c$-pairs of the form $([s_1, \delta_1]_{\tau_c}, [s_2, \delta_2]_{\tau_c})$ and $([s'_1, \delta_1]_{\tau_c}, [s'_2, \delta_2]_{\tau_c})$. Write $Z' = Z_G([s'_1, \delta_1]_{\tau_c}, [s'_2, \delta_2]_{\tau_c})$ which is a complex reductive group since the $c$-pair is reductive. Suppose that there is $g \in G$ such that
\[
([s_1, \delta_1]_{\tau_c}, [s_2, \delta_2]_{\tau_c}) = g ([s'_1, \delta_1]_{\tau_c}, [s'_2, \delta_2]_{\tau_c}) g^{-1}. 
\]
Then $S_c$ and $gS_cg^{-1}$ are Cartan subgroups of $Z'$, so there is an element $h \in Z'$ such that $h S_c h^{-1} = g S_c g^{-1}$ and then $g' = h^{-1}g$ is contained in $N_G(S_c)$. We have 
\[
g'([\id, \delta_1]_{\tau_c},[\id, \delta_2]_{\tau_c})(g')^{-1} = ([\id, \delta'_1]_{\tau_c}, [\id, \delta'_2]_{\tau_c}),
\]
where $(\delta'_1,\delta'_2)$ is an irreducible $c$-pair in $D_c$ and therefore, by Proposition \ref{co description of rank zero c-pairs for G}, there exists $\delta \in D_c$ such that $\delta (\delta'_1,\delta'_2) \delta^{-1} = (\delta_1,\delta_2)$. Noting that $[\id,\delta]_{\tau_c}$ conmutes with $S_c$ since $S_c$ is the centre of $Z_{G}(S_c)$, it follows that $g'' = [\id, \delta]_{\tau_c} \cdot g'\in N_G(S_c)$ and
\[
([s_1, \delta_1]_{\tau_c}, [s_2, \delta_2]_{\tau_c}) = ([g'' s'_1(g'')^{-1}, \delta_1]_{\tau_c}, [g''s'_2(g'')^{-1}, \delta_2]_{\tau_c}).
\]
Thus $(s_1,s_2)$ and $(s'_1,s'_2)$ define the same point of $(\ol{S}_c \times \ol{S}_c)/W_c$.
\qed

\bco
\label{co olzeta for G}
There is a bijective morphism
\beq
\label{eq ol zeta}
\quotient{\left( (\CC^* \times \CC^*) \otimes_\ZZ \Lambda_{\ol{S}_c} \right)}{W_c} \stackrel{1:1}{\lra} \Rr(G)_d.
\eeq 
and $\ol{\Rr(G)}_d=\quotient{\left((\CC^* \times \CC^*) \otimes_\ZZ \Lambda_{\ol{S}_c} \right)}{W_c}$ is the normalization of $\Rr(G)_d$.
\eco 

\bpf
Due to the isomorphism $\Theta_{\ol{S}_c} : \ol{S}_c \stackrel{\cong}{\lra} \CC^* \otimes_\ZZ \Lambda_{\ol{S}_c}$ defined in (\ref{eq definition of Theta}) and knowing that the action of $W_c$ commutes with it, the first statement follows  from Theorem \ref{tm homeomorphism RrG}.

The second statement follows from (\ref{eq ol zeta}) and Zariski's Main Theorem.
\epf

\brm
\label{rm Normality of RrG}
\nr{This is proved in \cite{thaddeus} for the case $d=0$. When the degree is trivial and $G = \GL(n,\CC)$ or $\SL(n,\CC)$, one obtains an isomorphism due to the normality of the target (see the discussion at the end of Section \ref{sc introduction}).}
\erm

\section{Hitchin equation and projectively flat bundles}
\label{sc connections}

Fix a maximal compact subgroup $K$ of $G$ and denote its Lie algebra by $\kK$. Take $\tau : \gG \to \gG$ to be the Cartan involution associated to the compact real form $\kK \subset \gG$. Then $\tau(k) = k$ and $\tau(ik) = -ik$ for every $k \in \kK$. 

Let $(E,\Phi)$ be a $G$-Higgs bundle and let $h$ be a metric on $E$, i.e. a $C^\infty$ reduction of $E$ to the maximal compact subgroup $K$ giving the $K$-bundle $E_h$. We define the involution on the adjoint bundle $\tau_h : E_h(\gG) \to E_h(\gG)$ using $\tau$ fibrewise.

Let $\dolbeault_E$ denote the Dolbeault operator of $E$ and set $A_h := \dolbeault_E + \tau_h(\dolbeault_E)$, which is the unique $K$-connection on $E_h$ compatible with $\dolbeault_E$, also known as the {\it Chern connection}. We denote by $F_h$ the curvature of $A_h$. 

Take the $C^\infty$ $(1,0)$-form $\dif x \in \Aa^{1,0}(X,\Oo)$ and $\dif \ol{x} \in \Aa^{0,1}(X,\Oo)$. Given a $G$-Higgs bundle $(E,\Phi)$, Hitchin introduced in \cite{hitchin-self_duality_equations} the following equation for a metric $h$ on $E$,
\beq
\label{eq Hitchin equations}
F_h + [\Phi \dif x, \tau_h(\Phi) \dif \ol{x}] = u \otimes \omega,
\eeq
where $u \in \zZ_{\gG}(\gG)$ and $\omega \in \Aa^2(X)$ is the volume form of the curve normalized to $2\pi i$. Recall that $u$ is determined by $d \in \pi_1(G)$.
 
In the elliptic case we have a splitting of the Hitchin equation.

\bpr
\label{co splitting of the hitchin equations}
If the $G$-Higgs bundle $(E,\Phi)$ is polystable then there exists a metric $h$ on $E$ that satisfies 
\[
F_h =  u \otimes \omega \qquad \nr{and} \qquad [\Phi \dif x, \tau_h(\Phi) \dif \ol{x}] = 0.
\]
\epr

\proof
By Corollary \ref{co E Phi polystable implies E polystable}, if the $G$-Higgs bundle $(E,\Phi)$ is polystable, then $E$ is polystable and by the Narasimhan--Seshadri--Ramanathan Theorem there exists a metric for which $F_h = u \otimes \omega$. 

By Corollary \ref{co description polystable G-Higgs bundles}, $(E,\Phi)$ is isomorphic to $(E^\rho,z \otimes \id_E)$ where $z \in \hH^{\omega_c}$. Then 
\[
[\Phi \dif x, \tau_h(\Phi) \dif \ol{x}] = [z,\tau(z)] \otimes \id_E \otimes (\dif x \wedge \dif \ol{x}) = 0
\]
since both $z$ and $\tau(z)$ belong to the abelian subalgebra $\hH$. 
\qed

One can easily show that a $G$-Higgs bundle $(E,\Phi)$ admitting a metric that satisfies (\ref{eq Hitchin equations}) is always polystable. Thus we see that Proposition \ref{co splitting of the hitchin equations} completes the proof of the Hitchin--Kobayashi correspondence in the elliptic case.

\bco
\label{co Hitchin Simpson correspondence}
A $G$-Higgs bundle $(E,\Phi)$ is polystable if and only if it admits a metric $h$ that satisfies the Hitchin equation (\ref{eq Hitchin equations}).
\eco

\brm
\label{rm Hitchin-Simpson easier}
\nr{Note that to prove the Hitchin--Kobayashi correspondence in the elliptic case we only make use of the Narasimhan--Seshadri--Ramanathan Theorem, the Jordan--Hölder reduction and Propositions \ref{pr E Phi semistable iff E semistable for G-Higgs bundles} and \ref{pr E Phi stable iff E stable for G-Higgs bundles}.} 
\erm

Let $\mathbf{E}_{G,d}$ be the (unique up to isomorphism) differentiable $G$-bundle of degree $d \in \pi_1(G)$ over the elliptic curve $X$. A $G$-connection $A$ on $\mathbf{E}_{G,d}$ is {\it flat} if the curvature vanishes, $F_A = 0$ (note that this forces $d = 0$). A $G$-connection $A$ on $\mathbf{E}_{G,d}$ is {\it projectively flat} or equivalently {\it $A$ has constant central curvature} if $F_A = a \otimes \omega$ for some $a \in \zZ_{\gG}(\gG)$. Due to topological considerations $a = u$, where $u \in \zZ_{\gG}(\gG)$ is determined by $d \in \pi_1(G)$. Let us denote by $\Cc(G)_d$ the moduli space of projectively flat connections on $\mathbf{E}_{G,d}$ and consequently $\Cc(G)_0$ is the moduli space of flat connections on $\mathbf{E}_{G,0}$.

We denote by $X^{\sharp}$ the moduli space of line bundles with flat connections over the elliptic curve $X$. Recalling that $T^*X \cong \Pic^0(X) \times H^0(X,\Omega^1_X)$, we have a homeomorpshism  
\beq
\label{eq homeomorphism between X sharp and T*X}
X^{\sharp} \stackrel{homeo}{\lra} T^*X 
\eeq
given by Hodge theory.

Let $S^\CC$ be a connected complex reductive abelian group. Recalling the isomorphism $\Theta_S$ given in (\ref{eq definition of Theta}), one can give a description of the moduli space of flat $S^\CC$-connections, denoted by $\Cc(S^\CC)_0$. Write $\mathbf{E}_{S,0}$ for the differentiable $S$-bundle with trivial topological class and recall that it is unique up to isomorphism.

Recall the isomorphism (\ref{eq definition of Theta}). For instance, the following result is contained in \cite[Theorem 9.10]{simpson2}.

\btm
Let $S^\CC$ be an abelian, connected complex Lie group. Then, the moduli space of flat $S^\CC$-connections over the elliptic curve $X$ is
\[
\Cc(S^\CC)_0 \cong X^\sharp \otimes_\ZZ \Lambda_S.
\]
\etm

Let $L \subset G$ be a reductive subgroup. We say that the $G$-connection $A$ {\it reduces} to the $L$-connection $A'$ when $A$ is gauge equivalent to the extension of structure group of $A'$ given by the natural injection $i : L \hookrightarrow G$.

Recall from (\ref{eq definition of q^pi}) that $d \in \pi_1(G)$ is determined by $(u,c) \in \pi_1(\ol{Z}) \times \pi_1(\ol{D})$, where $\pi_1(\ol{Z}) \subset \zZ_{\gG}(\gG)$ and $\pi_1(\ol{D}) = C$ as described in Section \ref{sc notation and results on Lie groups}. Take $L_c$ as defined in (\ref{eq definition of L_c}) and denote by $K_c$ its maximal compact subgroup.

\bpr
\label{co description of flat G-connections}
Every projectively flat connection $A$ on $\mathbf{E}_{G,d}$ reduces to a projectively flat $L_c$-connection. Futhermore $A$ is gauge equivalent to 
\[
A_{(\rho,z)} = A_\rho + z \dif x + \tau(z) \dif \ol{x}, 
\]
where $A_\rho$ is the Chern connection of $E^\rho$ given by $\rho : \Gamma_\RR \to K_c$ and $z \in \hH^{\omega_c}$.

The projectively flat connections $A_{(\rho,z)}$ and $A_{(\rho',z')}$ are gauge equivalent if and only if there exists $g \in K$ such that $(\rho',z') = (g \rho g^{-1}, \ad_g z)$.
\epr 

\proof
From a polystable $G$-Higgs bundle $(E,\Phi)$ we can construct a $G$-connection on $\mathbf{E}_{G,d}$ as follows
\[
A = A_h + \Phi \dif x + \tau_h(\Phi) \dif \ol{x}.
\]
Two isomorphic polystable $G$-Higgs bundles give rise to gauge equivalent flat $G$-connec\-tions. By Corollary \ref{co Hitchin Simpson correspondence}, the above $G$-connection is projectively flat if and only if $(E,\Phi)$ is po\-ly\-sta\-ble. The description of polystable $G$-Higgs bundles in Corollary \ref{co description polystable G-Higgs bundles} implies the proposition.
\qed

Denote by $\mathbf{E}_{L_c,\ell_d}$ the differentiable bundle underlying $E^{x_0}_{L_c,\ell_d}$, the $L_c$-bundle with degree $\ell_d$ defined in Remark \ref{rm definition of Ee_G}, and let $A^{x_0}_{L_c,\ell_d}$ be its Chern connection. Setting $p: X \times (X^\sharp \otimes \Lambda_{S_c}) \to X$, we define the family
\[
(\mathbf{F}'_{L_c,\ell_d} , (\Aa')_{L_c, \ell_d}) =\left( \mathbf{P}_{S,0} \otimes p^*\mathbf{E}_{L_c,\ell_d}, \Aa_{S_d,0} \otimes p^*A^{x_0}_{L_c, \ell_d} \right),
\]
noting that $S_c$ is the centre of $L_c$. This family is parametrized by $X^\sharp \otimes_\ZZ \Lambda_{S_c}$.

Recall $F_c$ and $\ol{S}_c$ as defined in \eqref{eq definition of olS_c and olD_c}. Let $J \in H^1(X,F_c)$ be a $F_c$-bundle and $A_J$ its Chern connection. By Corollary \ref{co H^1XF acts trivially on MstG}, one has the following.

\bpr
\label{co nabla_J otimes A is gauge equivalent to A}
Let $A$ be any $L_c$-connection on $\mathbf{E}_{L_c,\ell_d}$, then $A_J \otimes A$ is gauge equivalent to $A$.
\epr

As a consequence of Proposition \ref{co nabla_J otimes A is gauge equivalent to A}, it follows that $(\mathbf{F}'_{L_c,\ell_d}, (\Aa')_{L_c,\ell_d})$ induces a fa\-mi\-ly of $L_c$-connections parametrized by the quotient of $\Lambda_{S_c}$ by the subgroup associated to $H^1(X,F_c)$. This quotient is $\Lambda_{\ol{S}_c}$, and therefore we obtain a family parametrized by $X^{\sharp} \otimes_\ZZ \Lambda_{\ol{S}_c}$ that we denote by $(\mathbf{F}_{L_c,\ell_d},\Aa_{L_c,\ell_d})$.

Using the natural injection $i : L_c \hookrightarrow G$, we construct, by extension of structure group,
\[
(\mathbf{F}_{G,d},\Aa_{G,d}) = i_* (\mathbf{F}_{L_c,\ell_d},\Aa_{L_c,\ell_d}),
\]
a family of projectively flat $G$-connections, which is also parametrized by $X^\sharp \otimes_\ZZ \Lambda_{\ol{S}_c}$. 

\brm
\label{rm DdG is associated to HhG}
\nr{The flat $G$-connection parametrized by $(\mathbf{F}_{G,d},\Aa_{G,d})$ at the point $f \in X^\sharp \otimes_\ZZ \Lambda_{\ol{S}_c}$ is of the form $A_{(\rho,z)}$. It is therefore associated to the polystable $G$-Higgs bundle $(E^\rho,z \otimes 1_\Oo)$ parametrized by $\Hh_{G,d}$ at the point $(t,s) \in T^* X \otimes_\ZZ \Lambda_{\ol{S}_c}$, where $(t,s)$ is the image of $f$ under the homeomorphism (\ref{eq homeomorphism between X sharp and T*X}). Therefore, by Proposition \ref{pr eta/W_c is injective}, two points $f_1, f_2 \in X^\sharp \otimes_\ZZ \Lambda_{\ol{S}_c}$ parametrize gauge equivalent connections if $f_2 = \omega \cdot f_1$ for some $\omega \in W_c$.}
\erm

\btm
\label{tm DdG bijective Xsharp otimes Lambda}
There exists a bijective morphism
\beq
\label{eq tm DdG bijective Xsharp otimes Lambda}
\quotient{(X^\sharp \otimes_\ZZ \Lambda_{\ol{S}_c})}{W_c} \stackrel{1:1}{\lra} \Cc(G)_d
\eeq
and $\ol{\Cc(G)}_d=\quotient{\left(X^\sharp \otimes_\ZZ \Lambda_{\ol{S}_c} \right)}{W_c}$ is the normalization of $\Cc(G)_d$.
\etm 

\proof 
The family $(\mathbf{F}_{G,d},\Aa_{G,d})$ induces a morphism from the parametrizing space to the moduli space 
\[
X^\sharp \otimes_\ZZ \Lambda_{\ol{S}_c} \lra \Cc(G)_d, 
\]
which is surjective by Proposition \ref{co description of flat G-connections}. By Remark \ref{rm DdG is associated to HhG}, this surjection factors through (\ref{eq tm DdG bijective Xsharp otimes Lambda}) giving an injection.

The second statement follows from (\ref{eq tm DdG bijective Xsharp otimes Lambda}) and Zariski's Main Theorem.
\qed

\brm
\label{rm Normality of CcG}
\nr{This is proved in \cite{thaddeus} for the trivial degree case. In the case of $G = \GL(n,\CC)$ or $\SL(n,\CC)$ and $d = 0$, (\ref{eq tm DdG bijective Xsharp otimes Lambda}) is an isomorphism since the target is normal. Normality of $\Cc(G)_0$ follows from the Isosingularity Theorem \cite[Theorem 10.6]{simpson2} and normality of $\Rr(G)_0$.}
\erm

\end{document}